\documentclass[12pt,oneside,reqno]{amsart}
\usepackage{graphicx}
\usepackage{mathrsfs}
\usepackage{stmaryrd}
\usepackage{amsfonts}
\usepackage{enumerate,amsmath,amssymb,amsthm}
\usepackage{booktabs}
\usepackage{diagbox}
\usepackage{amsfonts,color}

\newcommand{\dif}{\mathrm{d}}

\newcommand{\be}{\begin{eqnarray}}
\newcommand{\ee}{\end{eqnarray}}
\newcommand{\ce}{\begin{eqnarray*}}
\newcommand{\de}{\end{eqnarray*}}
\newtheorem{theorem}{Theorem}[section]
\newtheorem{lemma}[theorem]{Lemma}
\newtheorem{remark}[theorem]{Remark}
\newtheorem{definition}[theorem]{Definition}
\newtheorem{proposition}[theorem]{Proposition}
\newtheorem{Examples}[theorem]{Examples}
\newtheorem{corollary}[theorem]{Corollary}
\def\e{\varepsilon}
\def\t{\theta}
\def\a{\alpha}

\def\b{\beta}

\def\p{\partial}
\def\g{\gamma}

\def\l{\lambda}
\def\la{\langle}
\def\ra{\rangle}
\def\[{{\Big[}}
\def\]{{\Big]}}
\def\<{{\langle}}
\def\>{{\rangle}}
\def\({{\Big(}}
\def\){{\Big)}}

\def\no{\nonumber}
\def\bt{\begin{theorem}}
\def\et{\end{theorem}}
\def\bl{\begin{lemma}}
\def\el{\end{lemma}}
\def\br{\begin{remark}}
\def\er{\end{remark}}
\def\bx{\begin{Examples}}
\def\ex{\end{Examples}}
\def\bd{\begin{definition}}
\def\ed{\end{definition}}
\def\bp{\begin{proposition}}
\def\ep{\end{proposition}}
\def\bc{\begin{corollary}}
\def\ec{\end{corollary}}

\def\cL{{\mathcal L}}

\def\cP{{\mathcal P}}

\def\mC{{\mathbb C}}

\def\mE{{\mathbb E}}

\def\mN{{\mathbb N}}

\def\mP{{\mathbb P}}

\def\mR{{\mathbb R}}

\def\mW{{\mathbb W}}

\def\sB{{\mathscr B}}
\def\sC{{\mathscr C}}

\def\sF{{\mathscr F}}

\def\sL{{\mathscr L}}

\def\geq{\geqslant}
\def\leq{\leqslant}

\begin{document}

\allowdisplaybreaks
\title{Asymptotic behaviors of multiscale McKean-Vlasov stochastic systems}

\author{Jie Xiang and Huijie Qiao}



\thanks{This work was supported by NSF of China (No.12071071) and the Jiangsu Provincial Scientific Research Center of Applied Mathematics (No. BK20233002).}

\thanks{Corresponding author: Huijie Qiao, hjqiaogean@seu.edu.cn}


\date{}

\dedicatory{Department of Mathematics,
Southeast University,\\
Nanjing, Jiangsu 211189, P.R.China\\
}

\begin{abstract}
In this paper, we investigate a class of multiscale McKean-Vlasov stochastic systems, where the entire system depends on the distributions of both fast and slow components. First of all, by applying the Poisson equation method, we prove that the slow component converges to the solution of the averaging equation in the $L^p$ ($p\geq 2$) space with the optimal convergence order $\frac12$. Then we establish a central limit theorem for these systems and derive the weak convergence rate using the Poisson equation technique and the regularity properties of the associated Cauchy problem.
\end{abstract}

\maketitle \rm

\section{Introduction}\label{intro}

McKean-Vlasov stochastic differential equations (SDEs for short), also referred to as distribution-dependent or mean-field SDEs, describe the evolution of individual particles within the mean field particle systems as the number of particles tends to infinity. A distinctive feature of these equations is that the coefficients depend not only on the solution process itself but also on its probability distribution. The study of McKean-Vlasov SDEs originated with McKean's foundational work \cite{hpm}, which was motivated by Kac's program in kinetic theory. Over the years, extensive investigations into McKean-Vlasov SDEs have generated notable advancements. These models have been analyzed from multiple perspectives, including well-posedness, stability, connections with nonlinear Fokker-Planck equations, exponential ergodicity, and more. We refer the interested reader to \cite{br, br1, dq1, dq2, hw, rz, Wang} and the references therein for a comprehensive overview.

Besides, multiscale stochastic systems, where the rates of change for different variables differ by orders of magnitude, are applied in various fields, such as chemistry, physics, climate dynamics and financial mathematics (See e.g. \cite{fp,mte,wtry}). For example, fast atmospheric and slow oceanic dynamics describe the climate evolution and state dynamics in electric power systems consist of fast- and slowly-varying elements.

Recently, multiscale McKean-Vlasov stochastic systems have garnered significant attention in the field of stochastic analysis. The inherent dependence on probability distributions, combined with the intricate interactions between slow and fast components, presents substantial challenges in analyzing the fundamental properties of these systems. Consequently, characterizing the behavior of such complex dynamics has emerged as a key research focus. A primary approach involves simplifying the original system by approximating it with a reduced model that captures its essential features. The averaging principle serves precisely this purpose by enabling effective dimension reduction and asymptotic approximation. In this paper, our first objective is to establish an averaging principle for a class of multiscale McKean-Vlasov stochastic systems in which the coefficients depend on the probability distributions of both the slow and fast components.

Concretely speaking, consider the following multiscale McKean-Vlasov stochastic system:
\be\left\{\begin{array}{l}
\dif X_t^{\e}=h_1(X_t^{\e}, \sL_{X_t^{\e}}, Y_t^{\e, y_0, \sL_{\xi}}, \sL_{Y_t^{\e, \xi}}) \dif t+\g_1(X_t^{\e}, \sL_{X_t^{\e}})\dif B_t, \\
X_0^{\e}=\varrho, \quad 0 \leq t\leq T, \\
\dif Y_t^{\e, \xi}=\frac{1}{\e} h_2(\sL_{X_t^{\e}}, Y_t^{\e, \xi}, \sL_{Y_t^{\e, \xi}}) \dif t+\frac{1}{\sqrt{\e}} \g_2(\sL_{X_t^{\e}}, Y_t^{\e,\xi}, \sL_{Y_t^{\e, \xi}})\dif W_t, \\
Y_0^{\e, \xi}=\xi, \quad 0 \leq t\leq T, \\
\dif Y_t^{\e, y_0, \sL_{\xi}}=\frac{1}{\e} h_2(\sL_{X_t^{\e}}, Y_t^{\e, y_0, \sL_{\xi}}, \sL_{Y_t^{\e, \xi}}) \dif t+\frac{1}{\sqrt{\e}} \g_2(\sL_{X_t^{\e}}, Y_t^{\e, y_0, \sL_{\xi}}, \sL_{Y_t^{\e, \xi}}) \dif W_t, \\
Y_0^{\e, y_0, \sL_{\xi}}=y_0, \quad 0 \leq t \leq T,
\end{array}\right.
\label{orieq}
\ee
where $\left(B_t\right),\left(W_t\right)$ are $d_1$-dimensional and $d_2$-dimensional standard Brownian motions, respectively, defined on the complete filtered probability space $(\Omega,\sF,\{\sF_t\}_{t \in[0, T]}, \mP)$ and are mutually independent. Moreover, these mappings
$h_1: \mR^n \times \cP_2\left(\mR^n\right) \times \mR^m \times \cP_2\left(\mR^m\right) \rightarrow \mR^n$,
$\g_1: \mR^n \times \cP_2\left(\mR^n\right) \rightarrow \mR^{n \times d_1}$,
$h_2: \cP_2\left(\mR^n\right)\times\mR^m \times \cP_2\left(\mR^m\right) \rightarrow \mR^m$,
$\g_2: \cP_2\left(\mR^n\right)\times\mR^m \times \cP_2\left(\mR^m\right) \rightarrow \mR^{m \times d_2}$
are all Borel measurable, and $\varrho, \xi$ are two random variables. Note that the system (\ref{orieq}) is more general than ones in some papers, such as \cite{bs, lwx, lx, qw, rsx, xq}. Here we demonstrate the following averaging principle for the system (\ref{orieq}) (See Theorem \ref{xbarxp})
\be
\mE\left(\sup_{0 \leq t \leq T}| X_t^{\e}- \bar{X}_t|^p\right)\leq C_T\e^{\frac p2}(1+\mE|\varrho|^{2p}+| y_0|^{2p}+\mE|\xi|^{2p}),
\label{xexbar}
\ee
where $X^\e, \bar X$ are solutions to the slow part of the system (\ref{orieq}) and the averaging equation \eqref{barxeq}, respectively. We mention that in \cite{lx}, for a general multiscale McKean-Vlasov stochastic system depending on the distribution of the fast component, Li and Xie obtained that (See \cite[Theorem 2.1]{lx})
\be
\sup_{0 \leq t \leq T}\mE|X_t^{\e}- \bar{X}_t|^p\leq C_T\e^{\frac p2}.\label{liave}
\ee
When the coefficients are independent of the distribution of the fast component, R\"{o}ckner, Sun and Xie \cite{rsx} established (\ref{liave}) for the case $p=2$. It is obvious that our result is stronger than those in \cite{lx, rsx, xq}.

Next, by (\ref{xexbar}), we know that $\{U^\e:=\frac{X^\e-\bar X}{\sqrt \e}, 0<\e<1\}$ is bounded in $L^p(\Omega,C([0,T],\mR^n))$. Therefore, it is natural to proceed with an investigation of the convergence for $U^\e$ as $\e\rightarrow 0$. This convergence result is called the central limit theorem (CLT). The foundational work on the CLT for multiscale SDEs is attributed to Khasminskii \cite{rk1}. Later, the CLT for multiscale SDEs has been greatly developed in \cite{rx, xie}. Additionally, \cite{ce, rxy, wr} explored the CLT in settings where the system transitions from finite to infinite dimensions. More recently, the CLT has gained significant attention in the context of multiscale McKean-Vlasov stochastic systems (cf. \cite{hlls, lx, xq}).

In this paper, our second objective is to establish a CLT for the system (\ref{orieq}). Specifically, we derive the limiting equation (\ref{uteq}) satisfied by $U$ and prove that $U^\e$ converges weakly to $U$, based on equations (\ref{ueeq}) and (\ref{uteq}). Furthermore, we establish a rate of weak convergence. A key challenge in obtaining this convergence rate arises from the fact that the coefficients of equation (\ref{uteq}) depend not only on probability distributions but also on derivatives with respect to these distributions. By employing the Poisson equation technique and leveraging the regularity properties of the associated Cauchy problem, we successfully overcome this difficulty.

It is worth noting that in \cite{xq}, we also demonstrated the CLT for a class of multiscale McKean-Vlasov SDEs, where the two fast components are independent of the distribution of the slow component. Nevertheless, in that study, we did not address the weak convergence rate, which constitutes one of the primary motivations for this work.

The remainder of the paper is organized as follows. Section \ref{noteassu} introduces the notations and assumptions used throughout the paper. The main results are presented in Section \ref{mare}. The proofs of the main results are provided in Section \ref{xbarxpproo} and \ref{cltthproo}, respectively. In Section \ref{example}, we give an example to illustrate the applicability of our results.

The following convention will be used throughout the paper: $C$ with or without indices will denote different positive constants whose values may change from one place to another.

\section{Notations and assumptions}\label{noteassu}

In this section, we will recall some notations and list all assumptions.

\subsection{Notations}\label{nn}
In this subsection, we introduce some notations used in the sequel.

We use $|\cdot|$ and $\|\cdot\|$ to denote the norms of vectors and matrices, respectively. Let $\<\cdot,\cdot\>$ represent the scalar product in $\mR^d$, and let $A^*$ denote the transpose of the matrix $A$.

Let $\sB(\mR^n)$ be the Borel $\sigma$-algebra on $\mR^n$ and $\cP(\mR^n)$ represent the space of all probability measures defined on $\sB(\mR^n)$ equipped with the usual topology of weak convergence. Let $\cP_2(\mR^n)$ be the collection of the probability measures $\mu$ on $\sB(\mR^n)$ satisfying
$$
\mu(|\cdot|^2):=\int_{\mR^n}|x|^2\mu(\dif x)<\infty.
$$
This space is a Polish space under the $L^2$-Wasserstein distance, defined as
$$
\mW_2(\mu_1,\mu_2):=\inf_{\pi\in\sC(\mu_1,\mu_2)}\bigg(\int_{\mR^n\times\mR^n}|x-y|^2\pi(\dif x,\dif y)\bigg)^{\frac{1}{2}},\mu_1,\mu_2 \in\cP_2(\mR^n),
$$
where $\sC(\mu_1,\mu_2)$ is the set of all couplings $\pi$ with marginals distributions $\mu_1$ and $\mu_2.$ Moreover, for any $x\in\mR^n$, the Dirac measure $\delta_x$ belongs to $\cP_2(\mR^n)$ and if $\mu_1=\sL_X,$ $\mu_2=\sL_Y$ are the corresponding distributions of random variables $X$ and $Y$ respectively, then
$$\mW_{2}(\mu_1,\mu_2)\leq(\mathbb{E}|X-Y|^{2})^{\frac12},$$
where $\mathbb{E}$ denotes the expectation with respect to $\mP$.

\subsection{Derivatives for functions on $\cP_2(\mR^d)$}\label{deri}
\quad In this subsection, we recall the definition of $L$-derivative for functions on $\cP_2(\mR^d)$. This definition was first introduced by Lions in \cite{card} (See also \cite{blpr,carm}), who employed abstract probability spaces to describe the $L$-derivatives. For the sake of clarity, we present a more straightforward formulation here (cf. \cite{rw}). Let $I$ be the identity map on $\mR^n$. For $\mu\in\cP_2(\mR^n)$ and $\phi\in L^2(\mR^n,\sB(\mR^n),\mu;\mR^n)$, where $L^2(\mR^n,\sB(\mR^n),\mu;\mR^n)$ stands for the space of Borel measurable functions $\phi: \mR^n\rightarrow\mR^n$ with $\int_{\mR^n}|\phi(x)|^2\mu(\dif x)<\infty$, $\mu(\phi):=\int_{\mR^n}\phi(x)\mu(\dif x)$. Moreover, simple calculations show that $\mu\circ(I+\phi)^{-1}\in\cP_2(\mR^n)$.

\bd\label{Lderidef}
(i) A function $g:\cP_2(\mR^n)\rightarrow\mR$ is called $L$-differentiable at $\mu\in\cP_2(\mR^n)$, if the functional
$$
L^2(\mR^n,\sB(\mR^n),\mu;\mR^n)\ni\phi\rightarrow g(\mu\circ(I+\phi)^{-1})
$$
is Fr\'{e}chet differentiable at $0\in L^2(\mR^n,\sB(\mR^n),\mu;\mR^n)$; that is, there exists a unique $\xi\in L^2(\mR^n,\sB(\mR^n),\mu;\mR^n)$ such that
$$
\lim_{\mu(|\phi|^2)\rightarrow0}\frac{g(\mu\circ(I+\phi)^{-1})-g(\mu)-\mu(\<\xi,\phi\>)}{\sqrt{\mu(|\phi|^2)}}=0.
$$
In this case, we denote $\p_\mu g(\mu)=\xi$ and call it the $L$-derivative of $g$ at $\mu$.\\
(ii) A function $g:\cP_2(\mR^n)\rightarrow\mR$ is called $L$-differentiable on $\cP_2(\mR^n)$ if $L$-derivative $\p_\mu g(\mu)$ exists for all $\mu\in\cP_2(\mR^n)$.
\ed

Let $d, l, k\in \mN=\{0,1,2,\ldots\}$. We introduce the following spaces of functions.
\begin{enumerate}[$\bullet$]
\item
The space $C^l_b(\mR^m,\mR)$.
A function $g(y)$ is said to be in $C^l_b(\mR^m,\mR)$, if $g(y)$ is $l$-times continuously differentiable and all its derivatives are bounded.

\item
The space $C^{(l,k)}_b(\cP_2(\mR^m),\mR)$.
A function $g(\nu)$ is said to be in $C^{(l,k)}_b(\cP_2(\mR^m),\mR)$, if $\nu\rightarrow g(\nu)$ is $l$-times continuously $L$-differentiable and all its $L$-derivatives are bounded. Additionally, we can find a version of $\p_\nu^lg(\nu)(\tilde{y}_1,\ldots,\tilde{y}_l)$ such that the mapping $(\tilde{y}_1,\ldots,\tilde{y}_l)\rightarrow\p_\nu^lg(\nu)(\tilde{y}_1,\ldots,\tilde{y}_l)$ is in $C^k_b(\mR^m\times\ldots\times\mR^m,\mR)$.

\item
The space $C^{2l,(l,l)}_b(\mR^m\times\cP_2(\mR^m),\mR)$.
A function $g(y,\nu)$ is said to be in $C^{2l,(l,l)}_b(\mR^m\times\cP_2(\mR^m),\mR)$, if for any $\nu\in\cP_2(\mR^m)$, the mapping $y\rightarrow g(y,\nu)$ is in $C^{2l}_b(\mR^m,\mR)$; for any $y\in\mR^m$, the mapping $\nu\rightarrow g(y,\nu)$ is in $C^{(l,l)}_b(\cP_2(\mR^m),\mR)$; and for any $1\leq i\leq l$, we can find a version of $\p_\nu^ig(y,\nu)(\tilde{y}_1,\ldots,\tilde{y}_i)$ such that the mapping $(y,\tilde{y}_1,\ldots,\tilde{y}_i)\rightarrow\p_\nu^ig(y,\nu)(\tilde{y}_1,\ldots,\tilde{y}_i)$ is in $C^{2l-i}_b(\mR^m\times\ldots\times\mR^m,\mR)$.

\item
The space $C^{(l,k),2l,(l,l)}_b(\cP_2(\mR^n)\times\mR^m\times\cP_2(\mR^m),\mR)$.
A function $g(\mu,y,\nu)$ is said to be in $C^{(l,k),2l,(l,l)}_b(\cP_2(\mR^n)\times\mR^m\times\cP_2(\mR^m),\mR)$, if for any $(y,\nu)\in\mR^m\times\cP_2(\mR^m)$, the mapping $\mu\rightarrow g(\mu,y,\nu)$ is in $C^{(l,k)}_b(\cP_2(\mR^n),\mR)$, and for any $\mu\in\cP_2(\mR^n)$, the mapping $(y,\nu)\rightarrow g(\mu,y,\nu)$ is in $C^{2l,(l,l)}_b(\mR^m\times\cP_2(\mR^m),\mR)$.

\item
The space $C^{d,(l,k),2l,(l,l)}_b(\mR^n\times\cP_2(\mR^n)\times\mR^m\times\cP_2(\mR^m),\mR)$.
A function $g(x,\mu,y,\nu)$ is said to be in $C^{d,(l,k),2l,(l,l)}_b(\mR^n\times\cP_2(\mR^n)\times\mR^m\times\cP_2(\mR^m),\mR)$, if for any $(\mu,y,\nu)\in\cP_2(\mR^n)\times\mR^m\times\cP_2(\mR^m)$, the mapping $x\rightarrow g(x,\mu,y,\nu)$ is in $C^d_b(\mR^n,\mR)$, and for any $x\in\mR^n$, the mapping $(\mu,y,\nu)\rightarrow g(x,\mu,y,\nu)$ is in $C^{(l,k),2l,(l,l)}_b(\cP_2(\mR^n)\times\mR^m\times\cP_2(\mR^m),\mR)$.

\item
The space
$\mC^{(l,k),2l,(l,l)}_b(\cP_2(\mR^n)\times\mR^m\times\cP_2(\mR^m),\mR)$.
A function $g(\mu,y,\nu)$ is said to be in $\mC^{(l,k),2l,(l,l)}_b(\cP_2(\mR^n)\times\mR^m\times\cP_2(\mR^m),\mR)$, if $g\in C^{(l,k),2l,(l,l)}_b(\cP_2(\mR^n)\times\mR^m\times\cP_2(\mR^m),\mR)$, and we can find a version of $\p_\mu^lg(\mu,y,\nu)(\tilde{x}_1,\ldots,\tilde{x}_l)$ such that the mapping $(y,\nu)\rightarrow\p_{\tilde{x}_1,\ldots,\tilde{x}_l}^k\p_\mu^lg(\mu,y,\nu)(\tilde{x}_1,\ldots,\tilde{x}_l)$ is in $C^{2l,(l,l)}_b(\mR^m\times\cP_2(\mR^m),\mR)$.

\item
The space $\mC^{d,(l,k),2l,(l,l)}_b(\mR^n\times\cP_2(\mR^n)\times\mR^m\times\cP_2(\mR^m),\mR)$.
A function $g(x,\mu,y,\nu)$ is said to be in $\mC^{d,(l,k),2l,(l,l)}_b(\mR^n\times\cP_2(\mR^n)\times\mR^m\times\cP_2(\mR^m),\mR)$, if $g\in C^{d,(l,k),2l,(l,l)}_b(\mR^n\times\cP_2(\mR^n)\times\mR^m\times\cP_2(\mR^m),\mR)$, and for any $x\in\mR^n$, the mapping $(\mu,y,\nu)\rightarrow g(x,\mu,y,\nu)$ is in $\mC^{(l,k),2l,(l,l)}_b(\cP_2(\mR^n)\times\mR^m\times\cP_2(\mR^m),\mR)$, and the mapping $(\mu,y,\nu)\rightarrow \p_x^dg(x,\mu,y,\nu)$ is in $\mC^{(l,k),2l,(l,l)}_b(\cP_2(\mR^n)\times\mR^m
\times\cP_2(\mR^m),\mR)$.

\item
The space $C_b^{1,2,(1,1),2,(1,1),2,(1,1)}(\mR_+\times\mR^n\times\cP_2(\mR^n)\times\mR^m\times\cP_2(\mR^m)\times\mR^n\times\cP_2(\mR^n),\mR)$.
A function $g(t,x,\mu,y,\nu,u,\pi)$ is said to be in $C_b^{1,2,(1,1),2,(1,1),2,(1,1)}(\mR_+\times\mR^n\times\cP_2(\mR^n)\times\mR^m\times\cP_2(\mR^m)\times\mR^n\times\cP_2(\mR^n),\mR)$, if for any $(x,\mu,y,\nu,u,\pi)\in \mR^n\times\cP_2(\mR^n)\times\mR^m\times\cP_2(\mR^m)\times\mR^n\times\cP_2(\mR^n)$, the mapping $t\rightarrow g(t,x,\mu,y,\nu,u,\pi)$ is in $C^{1}_b(\mR_+,\mR)$; for any $(t,u,\pi)\in \mR_+\times\mR^n\times\cP_2(\mR^n)$, the mapping $(x,\mu,y,\nu)\rightarrow g(t,x,\mu,y,\nu,u,\pi)$ is in $C^{2,(1,1),2,(1,1)}_b(\mR^n\times\cP_2(\mR^n)\times\mR^m\times\cP_2(\mR^m),\mR)$; for any $(t,x,\mu,y,\nu,\pi)\in\mR_+\times\mR^n\times\cP_2(\mR^n)
\times\mR^m\times\cP_2(\mR^m)\times\cP_2(\mR^n)$, the mapping $u\rightarrow g(t,x,\mu,y,\nu,u,\pi)$ is in $C^{2}_b(\mR^n,\mR)$; and for any $(t,x,\mu,y,\nu,u)\in\mR_+\times\mR^n\times\cP_2(\mR^n)
\times\mR^m\times\cP_2(\mR^m)\times\mR^n$, the mapping $\pi\rightarrow g(t,x,\mu,y,\nu,u,\pi)$ is in $C^{(1,1)}_b(\cP_2(\mR^n),\mR)$.
\end{enumerate}

\subsection{Assumptions}\label{ass}

In this subsection, we present all the assumptions used in the sequel:
\begin{enumerate}[$(\mathbf{H}^1_{h_{1}, \g_{1}})$]
\item
There exists a constant $L_{h_1, \g_1}>0$ such that for $x_i \in \mR^n$, $\mu_i \in \cP_2\left(\mR^n\right)$, $y_i \in \mR^m$, $\nu_i \in \cP_2\left(\mR^m\right), i=1,2$,
\ce
&&|h_1(x_1, \mu_1, y_1, \nu_1)-h_1(x_2, \mu_2, y_2, \nu_2)|^2+\|\g_1(x_1,\mu_1)-\g_1(x_2, \mu_2)\|^2 \\
&\leq& L_{h_1, \g_1}(|x_1-x_2|^2+\mW_2^2(\mu_1, \mu_2)+|y_1-y_2|^2+\mW_2^2(\nu_1, \nu_2)) .
\de
\end{enumerate}
\begin{enumerate}[$(\mathbf{H}^1_{h_{2}, \g_{2}})$]
\item
There exists a constant $L_{h_2, \g_2}>0$ such that for $\mu_i \in \cP_2\left(\mR^n\right)$, $y_i \in \mR^m$, $\nu_i \in \cP_2(\mR^m), i=1,2$,
\ce
&&|h_2(\mu_1, y_1, \nu_1)-h_2(\mu_2, y_2, \nu_2)|^2+\|\g_2(\mu_1, y_1, \nu_1)-\g_2(\mu_2, y_2, \nu_2)\|^2\\
&\leq& L_{h_2, \g_2}(\mW_2^2(\mu_1, \mu_2)+|y_1-y_2|^2+\mW_2^2(\nu_1, \nu_2)).
\de
\end{enumerate}
\begin{enumerate}[$(\mathbf{H}^2_{h_{2}, \g_{2}})$]
\item
For some $p\geq 1$, there exist two constants $\b_1>0, \b_2>0$ satisfying $\b_1-\b_2>4pL_{h_2, \g_2}$ such that for $\mu \in \cP_2\left(\mR^n\right)$, $y_i \in \mR^m$, $\nu_i \in \cP_2(\mR^m), i=1,2$,
\ce
&&2\< y_1-y_2, h_2(\mu, y_1, \nu_1)-h_2(\mu, y_2, \nu_2)\>+(2p-1)\|\g_2(\mu, y_1, \nu_1)-\g_2(\mu, y_2, \nu_2)\|^2\\
&\leq& -\b_1|y_1-y_2|^2+\b_2\mW_2^2(\nu_1, \nu_2).
\de
\end{enumerate}

\br\label{hgammarem}
$(i)$ $(\mathbf{H}^1_{h_{1}, \g_{1}})$ yields that
there exists a constant $\bar{L}_{h_1, \g_1}>0$ such that for $x \in \mR^n, \mu \in \cP_2\left(\mR^n\right), y \in \mR^m$, $\nu\in \cP_2\left(\mR^m\right)$,
$$
|h_1(x, \mu, y, \nu)|^2+\|\g_1(x, \mu)\|^2\leq \bar{L}_{h_1, \g_1}(1+|x|^2+\mu(|\cdot|^2)+|y|^2+\nu(|\cdot|^2)) .
$$

$(ii)$ $(\mathbf{H}_{h_2, \g_2}^1)$ implies that
there exists a constant $\bar{L}_{h_2, \g_2}>0$ such that for $\mu \in \cP_2\left(\mR^n\right)$, $y\in \mR^m$, $\nu\in \cP_2(\mR^m)$,
$$
|h_2(\mu, y, \nu)|^2+\|\g_2(\mu, y, \nu)\|^2\leq \bar{L}_{h_2, \g_2}(1+\mu(|\cdot|^2)+|y|^2+\nu(|\cdot|^2)) .
$$

$(iii)$ By $(\mathbf{H}_{h_2, \g_2}^1)$ and $(\mathbf{H}_{h_2, \g_2}^2)$, we can obtain
for $\mu \in \cP_2\left(\mR^n\right)$, $y \in \mR^m$, $\nu\in \cP_2\left(\mR^m\right)$,
$$
2\< y, h_2(\mu, y, \nu)\>+(2p-1)\|\g_2(\mu, y, \nu)\|^2\leq -\a_1|y|^2+\a_2\nu(|\cdot|^2)+C(1+\mu(|\cdot|^2)),
$$
where $\a_1:=\b_1-2pL_{h_2, \g_2}$, $\a_2:=\b_2+(2p-1)L_{h_2, \g_2}$, $\a_1-\a_2-L_{h_2, \g_2}>0$, and $C>0$ is a constant.

$(iv)$ $(\mathbf{H}_{h_2, \g_2}^2)$ assures the existence and uniqueness of invariant measures for the frozen equation. If we weaken this condition, invariant measures for the frozen equation probably are not unique and then there are more than one averaging equations. In the forthcoming work, we will consider this interesting case.
\er

\section{Main results}\label{mare}

In this section, we formulate the main result in this paper.

\subsection{The averaging principle for multiscale McKean-Vlasov SDEs}\label{avemare}

In this subsection, we present the averaging principle result for multiscale McKean-Vlasov SDEs.

Firstly, we recall the system \eqref{orieq}, i.e. for fixed $T>0$,
\ce\left\{\begin{array}{l}
\dif X_t^{\e}=h_1(X_t^{\e}, \sL_{X_t^{\e}}, Y_t^{\e, y_0, \sL_{\xi}}, \sL_{Y_t^{\e, \xi}}) \dif t+\g_1(X_t^{\e}, \sL_{X_t^{\e}}) \dif B_t, \\
X_0^{\e}=\varrho, \quad 0\leq t\leq T, \\
\dif Y_t^{\e, \xi}=\frac{1}{\e} h_2(\sL_{X_t^{\e}}, Y_t^{\e, \xi}, \sL_{Y_t^{\e, \xi}}) \dif t+\frac{1}{\sqrt{\e}} \g_2(\sL_{X_t^{\e}}, Y_t^{\e, \xi}, \sL_{Y_t^{\e, \xi}}) \dif W_t, \\
Y_0^{\e, \xi}=\xi, \quad 0\leq t\leq T, \\
\dif Y_t^{\e, y_0, \sL_{\xi}}=\frac{1}{\e} h_2(\sL_{X_t^{\e}}, Y_t^{\e, y_0, \sL_{\xi}}, \sL_{Y_t^{\e, \xi}}) \dif t+\frac{1}{\sqrt{\e}} \g_2(\sL_{X_t^{\e}}, Y_t^{\e, y_0, \sL_{\xi}}, \sL_{Y_t^{\e, \xi}}) \dif W_t, \\
Y_0^{\e, y_0, \sL_{\xi}}=y_0, \quad 0\leq t\leq T,
\end{array}
\right.
\de
where $\mE|\varrho|^{2p}<\infty$ and $\mE|\xi|^{2p}<\infty$ and $p$ is the same to that in $(\mathbf{H}_{h_2, \g_2}^2)$. Under ($\mathbf{H}_{h_1, \g_1}^1$) and ($\mathbf{H}_{h_2, \g_2}^1$),  it follows from \cite[Theorem 2.1]{Wang} that the system \eqref{orieq} admits a unique strong solution $(X_\cdot^{\e},Y_\cdot^{\e, \xi},Y_\cdot^{\e, y_0, \sL_{\xi}})$.

For fixed $\mu \in \cP_2\left(\mR^n\right)$, we consider the following system:
\be\left\{\begin{array}{l}
\dif Y_t^{\mu,\xi}=h_2(\mu,Y_t^{\mu,\xi}, \sL_{Y_t^{\mu,\xi}}) \dif t+\g_2(\mu,Y_t^{\mu,\xi}, \sL_{Y_t^{\mu,\xi}}) \dif W_t, \\
Y_0^{\mu,\xi}=\xi, \quad 0 \leq t \leq T, \\
\dif Y_t^{\mu,y_0, \sL_{\xi}}= h_2(\mu,Y_t^{\mu,y_0, \sL_{\xi}}, \sL_{Y_t^{\mu,\xi}})\dif t+\g_2(\mu,Y_t^{\mu,y_0, \sL_{\xi}}, \sL_{Y_t^{\mu,\xi}}) \dif W_t, \\
Y_0^{\mu,y_0, \sL_{\xi}}=y_0, \quad 0 \leq t \leq T.
\end{array}
\right.
\label{frozeq}
\ee
Under the assumption ($\mathbf{H}_{h_2, \g_2}^1$),  it is established in \cite[Theorem 2.1]{Wang} that the system \eqref{frozeq} admits a unique strong solution $(Y_\cdot^{\mu,\xi},Y_\cdot^{\mu, y_0, \sL_{\xi}})$. By analyzing the pair $(Y_\cdot^{\mu,y_0,\sL_{\xi}},\sL_{Y_\cdot^{\mu,\xi}})$, it follows from \cite[Theorem 4.11]{rrw} that one can construct a Markov process with the same distribution on a new probability space. Furthermore, under the assumption ($\mathbf{H}_{h_2, \g_2}^2$), \cite[Theorem 4.12]{rrw} ensures the existence of a unique invariant probability measure $\eta^\mu\times\delta_{\eta^\mu}$ for this Markov process, where $\eta^\mu$ is the unique invariant probability measure corresponding to the first equation in the system \eqref{frozeq} (See \cite[Theorem 3.1]{Wang}). This provides the basis for deriving the following averaging equation on the filtered probability space $(\Omega,\sF,\left\{\sF_t\right\}_{t \in[0, T]}, \mP)$:
\be\left\{\begin{array}{l}
\dif \bar{X}_t=\bar{h}_1(\bar{X}_t, \sL_{\bar{X}_t}) \dif t+\g_1(\bar{X}_t, \sL_{\bar{X}_t}) \dif B_t, \\
\bar{X}_0=\varrho,
\end{array}
\right.
\label{barxeq}
\ee
where $\bar{h}_1(x,\mu)=\int_{\mR^m\times\cP_2(\mR^m)}h_1(x,\mu,y,\nu)(\eta^\mu\times\delta_{\eta^\mu})(\dif y, \dif\nu)$. By the similar deduction to that for Theorem 3.1 in \cite{xq} or Lemma 4.7 in \cite{qw}, we know that Eq.\eqref{barxeq} has a unique solution $\bar{X}$.

The following result describes the relationship between $X^\e$ and $\bar X$, which is the first result in this paper.

\bt\label{xbarxp}
Supposed that $(\mathbf{H}_{h_1, \g_1}^1)$, $(\mathbf{H}_{h_2, \g_2}^1)$ and $(\mathbf{H}_{h_2, \g_2}^2)$ hold, $h_1\in C_b^{2,(1,1),2,(1,1)}(\mR^n\times\cP_2(\mR^n)\times\mR^m\times\cP_2(\mR^m),\mR^n)$, $h_2\in C_b^{(1,1),2,(1,1)}(\cP_2(\mR^n)\times\mR^m\times\cP_2(\mR^m),\mR^m)$ and $\g_2\in C_b^{(1,1),2,(1,1)}(\cP_2(\mR^n)\times\mR^m\times\cP_2(\mR^m),\mR^{m\times d_2})$. Then, for any $p\geq2$ $($p is the same as that in $(\mathbf{H}_{h_2, \g_2}^2)$$)$, there exists a positive constant $C_T$ independent of $\e$ such that
$$
\mE\left(\sup_{0 \leq t \leq T}| X_t^{\e}- \bar{X}_t|^p\right)\leq C_T\e^{\frac p2}(1+\mE|\varrho|^{2p}+| y_0|^{2p}+\mE|\xi|^{2p}).
$$
\et

The proof of Theorem \ref{xbarxp} is placed in Section \ref{xbarxpproo}.

\subsection{The central limit theorem for multiscale McKean-Vlasov SDEs}\label{cltmare}

In this subsection, we state the central limit theorem for multiscale McKean-Vlasov SDEs.

Now, we consider the deviation process $\{U_t^{\e}:=\frac{X_t^{\e}- \bar{X}_t}{\sqrt\e}\}_{t\geq0}$. By combining \eqref{orieq}, \eqref{barxeq} with Theorem \ref{xbarxp}, we know that $U^\e$ satisfies the following equation:
\be\left\{\begin{array}{l}
\dif U_t^{\e}=\frac{1}{\sqrt{\e}}[h_1(X_t^{\e}, \sL_{X_t^{\e}},  Y_t^{\e, y_0, \sL_{\xi}}, \sL_{Y_t^{\e, \xi}})-\bar{h}_1(\bar{X}_t, \sL_{\bar{X}_t})] \dif t \\
\qquad\quad+\frac{1}{\sqrt{\e}}[\g_1(X_t^{\e}, \sL_{X_t^{\e}})-\g_1(\bar{X}_t, \sL_{\bar{X}_t})] \dif B_t,\\
U_0^{\e}=0,
\end{array}
\right.
\label{ueeq}
\ee
and $\{U^\e, 0<\e<1\}$ is bounded in $L^p(\Omega,C([0,T],\mR^n))$.

Next, in order to study the limit of $U_t^{\e}$ as $\e\rightarrow 0$, we introduce some notations. Define the operator $\cL$ as follows: for every $\psi\in C_b^{2,(1,1)}(\mR^m\times\cP_2(\mR^m),\mR)$,
\ce
\cL\psi(y,\nu)
&:=&h_2(\mu,y,\nu)\cdot\p_y\psi(y,\nu)+\frac12Tr[\g_2\g_2^*(\mu,y,\nu)\cdot\p^2_{y}\psi(y,\nu)]\\
&&+\int_{\mR^m}\left[h_2(\mu,\tilde{y},\nu)\cdot\p_\nu\psi(y,\nu)(\tilde{y})+\frac12Tr[\g_2\g_2^*(\mu,\tilde{y},\nu)\cdot\p_{\tilde{y}}\p_{\nu}\psi(y,\nu)(\tilde{y})]\right]\nu(\dif\tilde{y}),
\de
and it is straightforward to verify that $\cL$ is the infinitesimal generator of $(Y_t^{\mu,y,\nu},\sL_{Y_t^{\mu,\xi}})$ with $\nu=\sL_{\xi}$. We now consider the following Poisson equation:
\be
-\cL\Psi(x,\mu,y,\nu)=h_1(x,\mu,y,\nu)-\bar{h}_1(x,\mu).
\label{h1poieq}
\ee
Under the assumptions of Theorem \ref{xbarxp} and based on \cite[Theorem 3.1]{lx}, Eq.\eqref{h1poieq} admits a unique solution $\Psi(x,\mu,y,\nu)$ in the space $C_b^{2,(1,1),2,(1,1)}(\mR^n\times\cP_2(\mR^n)\times\mR^m\times\cP_2(\mR^m),\mR^n)$. This solution is given by
\ce
\Psi(x,\mu,y,\nu)=\int_0^\infty\big[\mE h_1(x,\mu, Y_s^{\mu,y,\nu}, \sL_{Y_s^{\mu,\xi}})-\bar{h}_1(x,\mu)\big]\dif s, \quad\nu=\sL_{\xi}.
\de
Then we define
\ce
&&\p_y\Psi_{\g_2}(x,\mu,y,\nu):=\p_y\Psi(x,\mu,y,\nu)\cdot\g_2(\mu,y,\nu),\\
&&\overline{(\p_y\Psi_{\g_2})(\p_y\Psi_{\g_2})^*}(x,\mu):=\int_{\mR^m\times\cP_2(\mR^m)}(\p_y\Psi_{\g_2}(x,\mu,y,\nu))(\p_y\Psi_{\g_2}(x,\mu,y,\nu))^*(\eta^\mu\times\delta_{\eta^\mu})(\dif y, \dif\nu),\\
&&\Upsilon(x,\mu):=\left(\overline{(\p_y\Psi_{\g_2})(\p_y\Psi_{\g_2})^*}(x,\mu)\right)^{\frac12}.
\de
Based on the above notation, we construct the following McKean-Vlasov SDE:
\be\left\{\begin{array}{l}
\dif U_t=\p_x\bar{h}_1(\bar{X}_t, \sL_{\bar{X}_t})U_t\dif t+\tilde{\mE}[\p_\mu\bar{h}_1(\bar{X}_t, \sL_{\bar{X}_t})(\tilde{\bar{X}}_t)\tilde{U}_t]\dif t+\p_x\g_1(\bar{X}_t, \sL_{\bar{X}_t})U_t\dif B_t \\
\qquad\quad+\tilde{\mE}[\p_\mu\g_1(\bar{X}_t, \sL_{\bar{X}_t})(\tilde{\bar{X}}_t)\tilde{U}_t]\dif B_t+\Upsilon(\bar{X}_t, \sL_{\bar{X}_t})\dif V_t,\\
U_0=0,
\end{array}
\right.
\label{uteq}
\ee
where $V$ is a $n$-dimensional standard Brownian motion independent of $B$, the process $(\tilde{\bar{X}}_t,\tilde{U}_t)$ is a copy of $(\bar{X}_t,U_t)$ defined on a copy $(\tilde{\Omega},\tilde{\sF},\{\tilde{\sF}_t\}_{t\geq0},\tilde{\mP})$ of the original probability space $(\Omega,\sF,\{\sF_t\}_{t\geq0},\mP)$, and $\tilde{\mE}$ is the expectation taken with respect to $\tilde{\mP}$. Since Eq.(\ref{uteq}) is linear, there exists a unique solution $U$.

The following central limit theorem is the second result in this paper.

\bt\label{cltth}
Suppose that $(\mathbf{H}_{h_2, \g_2}^2)$ holds for $p\geq3$. Assume that $(\mathbf{H}_{h_1, \g_1}^1)$, $(\mathbf{H}_{h_2, \g_2}^1)$ hold and $\g_1\in \big(C_b^{4,(1,3)}\cap C_b^{4,(2,2)}\cap C_b^{4,(3,1)}\big)(\mR^n\times\cP_2(\mR^n),\mR^{n\times d_1})$, $h_2\in \big(C_b^{(3,1),6,(3,3)}\cap \mathbb{C}_b^{(1,3),4,(2,2)}\cap \mathbb{C}_b^{(2,2),2,(1,1)}\big)(\cP_2(\mR^n)\times\mR^m\times\cP_2(\mR^m),\mR^m), \g_2\in \big(C_b^{(3,1),6,(3,3)}\cap \mathbb{C}_b^{(1,3),4,(2,2)}\cap \mathbb{C}_b^{(2,2),2,(1,1)}\big)(\cP_2(\mR^n)\times\mR^m\times\cP_2(\mR^m),\mR^{m\times d_2})$ and
$h_1\in \big(C_b^{4,(3,1),6,(3,3)}\cap \mathbb{C}_b^{4,(1,3),4,(2,2)}\cap \mathbb{C}_b^{4,(2,2),2,(1,1)}\big)(\mR^n\times\cP_2(\mR^n)\times\mR^m\times\cP_2(\mR^m),\mR^n)$. Then the process $U^{\e}$ converges weakly to $U$ in $C([0,T];\mR^n)$ as $\e\rightarrow0$. Moreover, for any $\psi\in \big(C_b^{(1,3)}\cap C_b^{(2,2)}\cap C_b^{(3,1)}\big)(\cP_2(\mR^n),\mR)$, we have
\be
\sup_{0 \leq t \leq T}|\psi(\sL_{U_t^{\e}})-\psi(\sL_{U_t})|\leq C_T\sqrt\e.
\label{weakrate}
\ee
\et

The proof of Theorem \ref{cltth} is postponed to Section \ref{cltthproo}.

\section{Proof of Theorem \ref{xbarxp}}\label{xbarxpproo}

In this section, we prove Theorem \ref{xbarxp}. The first and second parts provide estimates for systems \eqref{orieq} and \eqref{frozeq}, respectively. In the third part, we conclude the proof of Theorem \ref{xbarxp} by integrating the results obtained in the preceding two subsections.

\subsection{Some estimates for the system \eqref{orieq}}\label{estiforxtyt}

In this subsection, we collect some estimates for the system \eqref{orieq}.

\bl\label{xtyt}
Suppose that assumptions $(\mathbf{H}_{h_1, \g_1}^1)$, $(\mathbf{H}_{h_2, \g_2}^1)$ and $(\mathbf{H}_{h_2, \g_2}^2)$ hold. Then, for any $p\geq1$ $($p is the same to that in $(\mathbf{H}_{h_2, \g_2}^2)$$)$, there exist positive constants $C$ and $C_T$ independent of $\e$ such that
\be
&&\mE\left(\sup_{0 \leq t \leq T}| X_t^{\e}|^{2p}\right)\leq C_T(1+\mE|\varrho|^{2p}+| y_0|^{2p}+\mE|\xi|^{2p}),\label{xeboun}\\
&&\sup_{0 \leq t \leq T}\mE| Y_t^{\e, \xi}|^{2p}\leq C_T(1+\mE|\varrho|^{2p}+| y_0|^{2p}+\mE|\xi|^{2p}),\label{yexibounw}\\
&&\sup_{0 \leq t \leq T}\mE| Y_t^{\e, y_0, \sL_{\xi}}|^{2p}\leq C_T(1+\mE|\varrho|^{2p}+| y_0|^{2p}+\mE|\xi|^{2p}),\label{yey0xiboun}\\
&&\mE\left(\sup_{0 \leq t \leq T}|Y_t^{\e, \xi}|^{2p}\right)\leq2\mE|\xi|^{2p}+\frac{C_T(1+\mE|\varrho|^{2p}+| y_0|^{2p}+\mE|\xi|^{2p})}{\e},\label{yexibounn}\\
&&\mE\left(\sup_{0 \leq t \leq T}|Y_t^{\e, y_0, \sL_{\xi}}|^{2p}\right)\leq2| y_0|^{2p}+\frac{C_T(1+\mE|\varrho|^{2p}+| y_0|^{2p}+\mE|\xi|^{2p})}{\e}. \label{yey0xibounn}
\ee
\el
\begin{proof}
By the same deduction in \cite[Lemma 4.1]{qw}, we can derive the estimates (\ref{xeboun})-(\ref{yey0xiboun}). In the following, we will establish the estimates (\ref{yexibounn}) and (\ref{yey0xibounn}).

For $Y_t^{\e, \xi}$, applying the It\^{o} formula to $|Y_t^{\e, \xi}|^{2p}$, we deduce that
\be
|Y_t^{\e, \xi}|^{2p}
&=&|\xi|^{2p}+\frac{2p}{\e}\int_0^t|Y_s^{\e, \xi}|^{2p-2}\< Y_s^{\e, \xi},h_2(\sL_{X_s^{\e}}, Y_s^{\e, \xi},\sL_{Y_s^{\e, \xi}})\>\dif s\no\\
&&+\frac{2p}{\sqrt{\e}}\int_0^t|Y_s^{\e, \xi}|^{2p-2}\< Y_s^{\e, \xi},\g_2(\sL_{X_s^{\e}},Y_s^{\e, \xi},\sL_{Y_s^{\e, \xi}})\dif W_s\>\no\\
&&+\frac{2p(p-1)}{\e}\int_0^t|Y_s^{\e, \xi}|^{2p-4}\|\g_2(\sL_{X_s^{\e}}, Y_s^{\e, \xi},\sL_{Y_s^{\e, \xi}})Y_s^{\e, \xi}\|^2\dif s\no\\
&&+\frac{p}{\e}\int_0^t|Y_s^{\e, \xi}|^{2p-2}\|\g_2(\sL_{X_s^{\e}}, Y_s^{\e, \xi},\sL_{Y_s^{\e, \xi}})\|^2\dif s.
\label{yexiito}
\ee
Note that \eqref{xeboun}, \eqref{yexibounw} and \eqref{yexiito} imply that
\ce
|Y_t^{\e, \xi}|^{2p}
&\leq&|\xi|^{2p}+\frac{p}{\e}\int_0^t|Y_s^{\e, \xi}|^{2p-2}\Big[2\< Y_s^{\e, \xi},h_2(\sL_{X_s^{\e}}, Y_s^{\e, \xi},\sL_{Y_s^{\e, \xi}})\>\\
&&\qquad\qquad\qquad\qquad\qquad\qquad+(2p-1)\|\g_2(\sL_{X_s^{\e}},Y_s^{\e, \xi},\sL_{Y_s^{\e, \xi}})\|^2\Big]\dif s\\
&&+\frac{2p}{\sqrt{\e}}\int_0^t|Y_s^{\e, \xi}|^{2p-2}\< Y_s^{\e, \xi},\g_2(\sL_{X_s^{\e}},Y_s^{\e, \xi},\sL_{Y_s^{\e, \xi}})\dif W_s\>\\
&\leq&|\xi|^{2p}+\frac{p}{\e}\int_0^t|Y_s^{\e, \xi}|^{2p-2}\left[-\a_1|Y_s^{\e, \xi}|^{2}+\a_2\mE|Y_s^{\e, \xi}|^{2}+C(1+\mE|X_s^{\e}|^{2}) \right]\dif s\\
&&+\frac{2p}{\sqrt{\e}}\int_0^t|Y_s^{\e, \xi}|^{2p-2}\< Y_s^{\e, \xi},\g_2(\sL_{X_s^{\e}},Y_s^{\e, \xi},\sL_{Y_s^{\e, \xi}})\dif W_s\>\\
&\leq&|\xi|^{2p}+\frac{p}{\e}\int_0^t\left[-(\a_1-\a_2-L_{h_2, \g_2})|Y_s^{\e, \xi}|^{2p}+C(1+\mE|X_s^{\e}|^{2p}+\mE|Y_s^{\e, \xi}|^{2p})\right]\dif s\\
&&+\frac{2p}{\sqrt{\e}}\int_0^t|Y_s^{\e, \xi}|^{2p-2}\< Y_s^{\e, \xi},\g_2(\sL_{X_s^{\e}},Y_s^{\e, \xi},\sL_{Y_s^{\e, \xi}})\dif W_s\>\\
&\leq&|\xi|^{2p}+\frac{C_T(1+\mE|\varrho|^{2p}+| y_0|^{2p}+\mE|\xi|^{2p})}{\e}\\
&&+\frac{2p}{\sqrt{\e}}\int_0^t|Y_s^{\e, \xi}|^{2p-2}\< Y_s^{\e, \xi},\g_2(\sL_{X_s^{\e}},Y_s^{\e, \xi},\sL_{Y_s^{\e, \xi}})\dif W_s\>.
\de
By the Burkholder-Davis-Gundy inequality, the Young inequality and Remark \ref{hgammarem}, we get
\ce
&&\mE\left(\sup_{0 \leq t \leq T}|Y_t^{\e, \xi}|^{2p}\right)\\
&\leq&\frac{C}{\sqrt{\e}}\mE\left[\int_0^T|Y_s^{\e, \xi}|^{4p-2}\|\g_2(\sL_{X_s^{\e}},Y_s^{\e, \xi},\sL_{Y_s^{\e, \xi}})\|^2\dif s\right]^{\frac12}\\
&&+\mE|\xi|^{2p}+\frac{C_T(1+\mE|\varrho|^{2p}+| y_0|^{2p}+\mE|\xi|^{2p})}{\e}\\
&\leq&\frac{C}{\sqrt{\e}}\mE\left[\sup_{0 \leq t \leq T}|Y_t^{\e, \xi}|^{2p}\int_0^T|Y_s^{\e, \xi}|^{2p-2}(1+\mE|X_s^{\e}|^{2}+|Y_s^{\e, \xi}|^{2}+\mE|Y_s^{\e, \xi}|^{2})\dif s\right]^{\frac12}\\
&&+\mE|\xi|^{2p}+\frac{C_T(1+\mE|\varrho|^{2p}+| y_0|^{2p}+\mE|\xi|^{2p})}{\e}\\
&\leq&\frac12\mE\left(\sup_{0 \leq t \leq T}|Y_t^{\e, \xi}|^{2p}\right)+\frac{C}{\e}\int_0^T(1+\mE|X_s^{\e}|^{2p}+\mE|Y_s^{\e, \xi}|^{2p})\dif s\\
&&+\mE|\xi|^{2p}+\frac{C_T(1+\mE|\varrho|^{2p}+| y_0|^{2p}+\mE|\xi|^{2p})}{\e}\\
&\leq&\frac12\mE\left(\sup_{0 \leq t \leq T}|Y_t^{\e, \xi}|^{2p}\right)+\mE|\xi|^{2p}+\frac{C_T(1+\mE|\varrho|^{2p}+| y_0|^{2p}+\mE|\xi|^{2p})}{\e},
\de
which implies that
\ce
\mE\left(\sup_{0 \leq t \leq T}|Y_t^{\e, \xi}|^{2p}\right)\leq2\mE|\xi|^{2p}+\frac{C_T(1+\mE|\varrho|^{2p}+| y_0|^{2p}+\mE|\xi|^{2p})}{\e}.
\de

Moreover, for $Y_t^{\e, y_0, \sL_{\xi}}$, using a similar argument as for $Y_t^{\e, \xi}$, we obtain the estimate (\ref{yey0xibounn}). The proof is complete.
\end{proof}

\subsection{Some estimates for the system \eqref{frozeq}}\label{estiforymu}

In this subsection, we present some estimates for the system \eqref{frozeq}.

By arguments similar to that in \cite[Lemma 4.6]{qw}, we obtain the following result.

\bl\label{meh1}
Suppose that assumptions $(\mathbf{H}_{h_1, \g_1}^1)$, $(\mathbf{H}_{h_2, \g_2}^1)$ and $(\mathbf{H}_{h_2, \g_2}^2)$ hold. Then it holds that for any $t\geq0$, $x\in \mR^n$, $\mu\in\cP_2\left(\mR^n\right)$, there exists a constant $C>0$ such that
\ce
&&|\mE h_1(x,\mu, Y_t^{\mu, y_0, \sL_{\xi}}, \sL_{Y_t^{\mu, \xi}})-\bar{h}_1(x,\mu)|^2\\
&\leq& Ce^{-(\b_1-\b_2-L_{h_2, \g_2})t}\left(1+\mu(|\cdot|^2)+|y_0|^2+\sL_{\xi}(|\cdot|^2)\right).
\de
\el

\subsection{The proof of Theorem \ref{xbarxp}}\label{xbarxpproode}
First of all, we prove the Lipschitz continuity of the coefficient $\bar{h}_1$.
By Lemma \ref{meh1} and ($\mathbf{H}_{h_1, \g_1}^1$), for any $x_i\in \mR^n$, $\mu_i\in\cP_2\left(\mR^n\right)$, $i=1,2$, we have
\ce
&&|\bar{h}_1(x_1,\mu_1)-\bar{h}_1(x_2,\mu_2)|^2\\
&\leq&3|\bar{h}_1(x_1,\mu_1)-\mE h_1(x_1,\mu_1, Y_t^{\mu_1, y_0, \sL_{\xi}}, \sL_{Y_t^{\mu_1, \xi}})|^2\\
&&+3|\mE h_1(x_2,\mu_2, Y_t^{\mu_2, y_0, \sL_{\xi}}, \sL_{Y_t^{\mu_2, \xi}})-\bar{h}_1(x_2,\mu_2)|^2\\
&&+3|\mE h_1(x_1,\mu_1, Y_t^{\mu_1, y_0, \sL_{\xi}}, \sL_{Y_t^{\mu_1, \xi}})-\mE h_1(x_2,\mu_2, Y_t^{\mu_2, y_0, \sL_{\xi}}, \sL_{Y_t^{\mu_2, \xi}})|^2\\
&\leq&Ce^{-(\b_1-\b_2-L_{h_2, \g_2})t}\left(1+\mu_1(|\cdot|^2)+\mu_2(|\cdot|^2)+|y_0|^2+\sL_{\xi}(|\cdot|^2)\right)\\
&&+C(|x_1-x_2|^2+\mW_2^2(\mu_1, \mu_2)).
\de
Letting $t\rightarrow\infty$, it immediately leads to
\be
|\bar{h}_1(x_1,\mu_1)-\bar{h}_1(x_2,\mu_2)|^2\leq C(|x_1-x_2|^2+\mW_2^2(\mu_1, \mu_2)).\label{barh1lip}
\ee

Note that
\ce
X_t^{\e}-\bar{X}_t
&=&\int_0^t\Big(h_1(X_s^{\e}, \sL_{X_s^{\e}}, Y_s^{\e, y_0, \sL_{\xi}}, \sL_{Y_s^{\e, \xi}})-\bar{h}_1(\bar{X}_s, \sL_{\bar{X}_s})\Big)\dif s\\
&&+\int_0^t\Big(\g_1(X_s^{\e}, \sL_{X_s^{\e}})-\g_1(\bar{X}_s, \sL_{\bar{X}_s})\Big) \dif B_s\\
&=&\int_0^t\Big(h_1(X_s^{\e}, \sL_{X_s^{\e}},Y_s^{\e, y_0, \sL_{\xi}}, \sL_{Y_s^{\e, \xi}})-\bar{h}_1(X_s^{\e}, \sL_{X_s^{\e}})\Big)\dif s\\
&&+\int_0^t\Big(\bar{h}_1(X_s^{\e}, \sL_{X_s^{\e}})-\bar{h}_1(\bar{X}_s, \sL_{\bar{X}_s})\Big)\dif s\\
&&+\int_0^t\Big(\g_1(X_s^{\e}, \sL_{X_s^{\e}})-\g_1(\bar{X}_s, \sL_{\bar{X}_s})\Big) \dif B_s.
\de
From the H\"{o}lder inequality, the Burkholder-Davis-Gundy inequality, ($\mathbf{H}_{h_1, \g_1}^1$) and \eqref{barh1lip}, it follows that for any $p\geq2$,
\ce
&&\mE\left(\sup_{0\leq t\leq T}|X_t^{\e}-\bar{X}_t|^p\right)\\
&\leq&C\mE\left(\sup_{0\leq t\leq T}\left|\int_0^t\Big(h_1(X_s^{\e}, \sL_{X_s^{\e}},Y_s^{\e, y_0, \sL_{\xi}}, \sL_{Y_s^{\e, \xi}})-\bar{h}_1(X_s^{\e}, \sL_{X_s^{\e}})\Big)\dif s\right|^p\right)\\
&&+CT^{p-1}\mE\int_0^T\left|\bar{h}_1(X_s^{\e}, \sL_{X_s^{\e}})-\bar{h}_1(\bar{X}_s, \sL_{\bar{X}_s})\right|^p\dif s\\
&&+C\mE\left[\int_0^T\left\|\g_1(X_s^{\e}, \sL_{X_s^{\e}})-\g_1(\bar{X}_s, \sL_{\bar{X}_s})\right\|^2 \dif s\right]^{\frac p2}\\
&\leq&C\mE\left(\sup_{0\leq t\leq T}\left|\int_0^t\Big(h_1(X_s^{\e}, \sL_{X_s^{\e}},Y_s^{\e, y_0, \sL_{\xi}}, \sL_{Y_s^{\e, \xi}})-\bar{h}_1(X_s^{\e}, \sL_{X_s^{\e}})\Big)\dif s\right|^p\right)\\
&&+C_T\int_0^T\mE|X_s^{\e}-\bar{X}_s|^p\dif s.
\de
The Gronwall inequality yields that
\be
&&\mE\left(\sup_{0\leq t\leq T}|X_t^{\e}-\bar{X}_t|^p\right)\no\\
&\leq&C_T\mE\left(\sup_{0\leq t\leq T}\left|\int_0^t\Big(h_1(X_s^{\e}, \sL_{X_s^{\e}},Y_s^{\e, y_0, \sL_{\xi}}, \sL_{Y_s^{\e, \xi}})-\bar{h}_1(X_s^{\e}, \sL_{X_s^{\e}})\Big)\dif s\right|^p\right).
\label{xebarxgrow}
\ee

In the following, we aim to estimate the right-hand part of the inequality above.
Applying the It\^{o} formula to
$\Psi(X_t^{\e}, \sL_{X_t^{\e}},Y_t^{\e, y_0, \sL_{\xi}}, \sL_{Y_t^{\e, \xi}})$, we have
\ce
&&\Psi(X_t^{\e}, \sL_{X_t^{\e}},Y_t^{\e, y_0, \sL_{\xi}}, \sL_{Y_t^{\e, \xi}})\\
&=&\Psi(\varrho,\sL_{\varrho},y_0,\sL_{\xi})+\int_0^th_1(X_s^{\e}, \sL_{X_s^{\e}},Y_s^{\e, y_0, \sL_{\xi}}, \sL_{Y_s^{\e, \xi}})\cdot\p_x\Psi(X_s^{\e}, \sL_{X_s^{\e}},Y_s^{\e, y_0, \sL_{\xi}}, \sL_{Y_s^{\e, \xi}})\dif s\no\\
&&+\frac12\int_0^tTr\left[\g_1\g_1^*(X_s^{\e}, \sL_{X_s^{\e}})\cdot\p_{x}^2\Psi(X_s^{\e}, \sL_{X_s^{\e}},Y_s^{\e, y_0, \sL_{\xi}}, \sL_{Y_s^{\e, \xi}})\right]\dif s\no\\
&&+\int_0^t\p_x\Psi(X_s^{\e}, \sL_{X_s^{\e}},Y_s^{\e, y_0, \sL_{\xi}}, \sL_{Y_s^{\e, \xi}})\cdot\g_1(X_s^{\e}, \sL_{X_s^{\e}})\dif B_s\no\\
&&+\frac1{\e}\int_0^t h_2(\sL_{X_s^{\e}}, Y_s^{\e, y_0, \sL_{\xi}},\sL_{Y_s^{\e, \xi}})\cdot\p_y\Psi(X_s^{\e}, \sL_{X_s^{\e}},Y_s^{\e, y_0, \sL_{\xi}}, \sL_{Y_s^{\e, \xi}})\dif s\no\\
&&+\frac1{2\e}\int_0^tTr\left[\g_2\g_2^*(\sL_{X_s^{\e}}, Y_s^{\e, y_0, \sL_{\xi}}, \sL_{Y_s^{\e, \xi}})\cdot\p_{y}^2\Psi(X_s^{\e}, \sL_{X_s^{\e}},Y_s^{\e, y_0, \sL_{\xi}}, \sL_{Y_s^{\e, \xi}})\right]\dif s\no\\
&&+\frac1{\sqrt\e}\int_0^t\p_y\Psi(X_s^{\e}, \sL_{X_s^{\e}},Y_s^{\e, y_0, \sL_{\xi}}, \sL_{Y_s^{\e, \xi}})\cdot\g_2(\sL_{X_s^{\e}}, Y_s^{\e, y_0, \sL_{\xi}}, \sL_{Y_s^{\e, \xi}})\dif W_s\no\\
&&+\mathbb{\tilde{E}}\int_0^t h_1(\tilde{X}_s^{\e}, \sL_{X_s^{\e}},\tilde{Y}_s^{\e, y_0, \sL_{\xi}}, \sL_{Y_s^{\e, \xi}})\cdot\p_\mu\Psi(X_s^{\e}, \sL_{X_s^{\e}},Y_s^{\e, y_0, \sL_{\xi}}, \sL_{Y_s^{\e, \xi}})(\tilde{X}_s^{\e})  \dif s\no\\
&&+\frac12\mathbb{\tilde{E}}\int_0^tTr\left[\g_1\g_1^*(\tilde{X}_s^{\e}, \sL_{X_s^{\e}})\cdot\p_{\tilde{x}}\p_\mu\Psi(X_s^{\e}, \sL_{X_s^{\e}},Y_s^{\e, y_0, \sL_{\xi}}, \sL_{Y_s^{\e, \xi}})(\tilde{X}_s^{\e})\right]\dif s\no\\
&&+\frac1\e\mathbb{\tilde{E}}\int_0^th_2(\sL_{X_s^{\e}}, \tilde{Y}_s^{\e, \tilde{\xi}}, \sL_{Y_s^{\e,\xi}})\cdot\p_\nu\Psi(X_s^{\e}, \sL_{X_s^{\e}},Y_s^{\e, y_0, \sL_{\xi}}, \sL_{Y_s^{\e, \xi}})(\tilde{Y}_s^{\e, \tilde{\xi}})\dif s\no\\
&&+\frac1{2\e}\mathbb{\tilde{E}}\int_0^tTr\left[\g_2\g_2^*(\sL_{X_s^{\e}}, \tilde{Y}_s^{\e, \tilde{\xi}}, \sL_{Y_s^{\e,\xi}})\cdot\p_{\tilde{y}}\p_\nu\Psi(X_s^{\e}, \sL_{X_s^{\e}},Y_s^{\e, y_0, \sL_{\xi}}, \sL_{Y_s^{\e, \xi}})(\tilde{Y}_s^{\e, \tilde{\xi}})\right]\dif s\no\\
&=&\Psi(\varrho,\sL_{\varrho},y_0,\sL_{\xi})+M_t^{\e,1}+\frac1{\sqrt\e}M_t^{\e,2}\no\\
&&+\mathbb{\tilde{E}}\int_0^t h_1(\tilde{X}_s^{\e}, \sL_{X_s^{\e}},\tilde{Y}_s^{\e, y_0, \sL_{\xi}}, \sL_{Y_s^{\e, \xi}})\cdot\p_\mu\Psi(X_s^{\e}, \sL_{X_s^{\e}},Y_s^{\e, y_0, \sL_{\xi}}, \sL_{Y_s^{\e, \xi}})(\tilde{X}_s^{\e})  \dif s\no\\
&&+\frac12\mathbb{\tilde{E}}\int_0^tTr\left[\g_1\g_1^*(\tilde{X}_s^{\e}, \sL_{X_s^{\e}})\cdot\p_{\tilde{x}}\p_\mu\Psi(X_s^{\e}, \sL_{X_s^{\e}},Y_s^{\e, y_0, \sL_{\xi}}, \sL_{Y_s^{\e, \xi}})(\tilde{X}_s^{\e})\right]\dif s\no\\
&&+\int_0^t\cL_1\Psi(X_s^{\e}, \sL_{X_s^{\e}},Y_s^{\e, y_0, \sL_{\xi}}, \sL_{Y_s^{\e, \xi}})\dif s+\frac1{\e}\int_0^t\cL\Psi(X_s^{\e}, \sL_{X_s^{\e}},Y_s^{\e, y_0, \sL_{\xi}}, \sL_{Y_s^{\e, \xi}})\dif s,\no\\
\de
where the process $(\tilde{X}_s^{\e},\tilde{Y}_s^{\e, y_0, \sL_{\xi}},\tilde{Y}_s^{\e, \tilde{\xi}})$ represents a copy of the original process $(X_s^{\e},Y_s^{\e, y_0, \sL_{\xi}},$ $Y_s^{\e, \xi})$, constructed on a probability space $(\tilde{\Omega},\tilde{\sF},\{\tilde{\sF}_t\}_{t \in[0, T]}, \tilde{\mP})$, which is an exact copy of the original probability space $(\Omega,\sF,\{\sF_t\}_{t \in[0, T]}, \mP)$, $\tilde{\mE}$ is the expectation taken with respect to $\tilde{\mP}$,
\ce
&&M_t^{\e,1}:=\int_0^t\p_x\Psi(X_s^{\e}, \sL_{X_s^{\e}},Y_s^{\e, y_0, \sL_{\xi}}, \sL_{Y_s^{\e, \xi}})\cdot\g_1(X_s^{\e}, \sL_{X_s^{\e}})\dif B_s,\\
&&M_t^{\e,2}:=\int_0^t\p_y\Psi(X_s^{\e}, \sL_{X_s^{\e}},Y_s^{\e, y_0, \sL_{\xi}}, \sL_{Y_s^{\e, \xi}})\cdot\g_2(\sL_{X_s^{\e}}, Y_s^{\e, y_0, \sL_{\xi}}, \sL_{Y_s^{\e, \xi}})\dif W_s,
\de
and
\ce
\cL_1:=\sum_{i=1}^nh_{1,i}(x,\mu,y,\nu)\frac{\p}{\p_{x_i}}+\frac12\sum_{i,j=1}^n\big(\g_1\g_1^*(x,\mu)\big)_{ij}\frac{\p^2}{\p_{x_i}\p_{x_j}}.
\de
Then based on the fact that $\Psi(x,\mu,y,\nu)$ is the unique solution to Eq.\eqref{h1poieq}, we can deduce that
\ce
&&\mE\left(\sup_{0\leq t\leq T}\left|\int_0^t\Big(h_1(X_s^{\e}, \sL_{X_s^{\e}},Y_s^{\e, y_0, \sL_{\xi}}, \sL_{Y_s^{\e, \xi}})-\bar{h}_1(X_s^{\e}, \sL_{X_s^{\e}})\Big)\dif s\right|^p\right)\\
&=&\mE\left(\sup_{0\leq t\leq T}\left|\int_0^t\cL\Psi(X_s^{\e}, \sL_{X_s^{\e}},Y_s^{\e, y_0, \sL_{\xi}}, \sL_{Y_s^{\e, \xi}})\dif s\right|^p\right)\\
&\leq&C\e^p\mE\left(\sup_{0\leq t\leq T}\Big|\Psi(X_t^{\e}, \sL_{X_t^{\e}},Y_t^{\e, y_0, \sL_{\xi}}, \sL_{Y_t^{\e, \xi}})\Big|^p\right)+C\e^p\mE\left|\Psi(\varrho,\sL_{\varrho},y_0,\sL_{\xi})\right|^p\\
&&+C\e^p\mE\mathbb{\tilde{E}}\Bigg(\sup_{0\leq t\leq T}\Big|\int_0^t h_1(\tilde{X}_s^{\e}, \sL_{X_s^{\e}},\tilde{Y}_s^{\e, y_0, \sL_{\xi}}, \sL_{Y_s^{\e, \xi}}),\\
&&\qquad\qquad\qquad\qquad\qquad\qquad\cdot\p_\mu\Psi(X_s^{\e}, \sL_{X_s^{\e}},Y_s^{\e, y_0, \sL_{\xi}}, \sL_{Y_s^{\e, \xi}})(\tilde{X}_s^{\e}) \dif s\Big|^p\Bigg)\\
&&+C\e^p\mE\mathbb{\tilde{E}}\Bigg(\sup_{0\leq t\leq T}\Big|\int_0^tTr\Big[\g_1\g_1^*(\tilde{X}_s^{\e}, \sL_{X_s^{\e}})\\
&&\qquad\qquad\qquad\qquad\qquad\qquad\cdot\p_{\tilde{x}}\p_\mu\Psi(X_s^{\e}, \sL_{X_s^{\e}},Y_s^{\e, y_0, \sL_{\xi}}, \sL_{Y_s^{\e, \xi}})(\tilde{X}_s^{\e})\Big]\dif s\Big|^p\Bigg)\\
&&+C\e^p\mE\left(\sup_{0\leq t\leq T}\Big|\int_0^t\cL_1\Psi(X_s^{\e}, \sL_{X_s^{\e}},Y_s^{\e, y_0, \sL_{\xi}}, \sL_{Y_s^{\e, \xi}})\dif s\Big|^p\right)\\
&&+C\e^p\mE\left(\sup_{0\leq t\leq T}|M_t^{\e,1}|^p\right)+C\e^{\frac p2}\mE\left(\sup_{0\leq t\leq T}|M_t^{\e,2}|^p\right)\\
&=&\sum_{i=1}^7\mathscr{J}_i.
\de

For $\mathscr{J}_1$ and $\mathscr{J}_2$, Lemma \ref{meh1} implies that
\ce
|\Psi(x,\mu,y,\nu)|
&\leq&\int_0^\infty\Big|\mE h_1(x,\mu, Y_s^{\mu,y,\nu}, \sL_{Y_s^{\mu,\xi}})-\bar{h}_1(x,\mu)\Big|\dif s\no\\
&\leq&\int_0^\infty Ce^{-\frac{\b_1-\b_2-L_{h_2, \g_2}}2s}\left(1+\mu(|\cdot|^2)+|y|^2+\nu(|\cdot|^2)\right)^{1/2}\dif s\no\\
&\leq&C\left(1+\mu(|\cdot|^2)+|y|^2+\nu(|\cdot|^2)\right)^{1/2},
\de
where $\nu=\sL_{\xi}$. Combining this with Lemma \ref{xtyt}, one can obtain that
\ce
\mathscr{J}_1
&\leq& C\e^p\left\{1+\mE\left(\sup_{0\leq t\leq T}|X_t^{\e}|^p\right)+\mE\left(\sup_{0\leq t\leq T}|Y_t^{\e, y_0, \sL_{\xi}}|^p\right)+\mE\left(\sup_{0\leq t\leq T}|Y_t^{\e, \xi}|^p\right)\right\}\\
&\leq& C_T\e^{p-1}\{1+\mE|\varrho|^{2p}+|y_0|^{2p}+\mE|\xi|^{2p}\},
\de
and
\ce
\mathscr{J}_2\leq C\e^p\{1+\mE|\varrho|^{p}+|y_0|^p+\mE|\xi|^p\}.
\de

For $\mathscr{J}_3$ and $\mathscr{J}_4$, by applying the H\"{o}lder inequality and Remark \ref{hgammarem}, we obtain
\ce
\mathscr{J}_3&\leq&C\e^pT^{p-1}\mE\mathbb{\tilde{E}}\Big[\int_0^T|h_1(\tilde{X}_s^{\e}, \sL_{X_s^{\e}},\tilde{Y}_s^{\e, y_0, \sL_{\xi}}, \sL_{Y_s^{\e, \xi}})|^p \dif s\Big]\\
&\leq&C\e^pT^p\left\{1+\mE\left(\sup_{0\leq t\leq T}|X_t^{\e}|^{p}\right)+\sup_{0\leq t\leq T}\mE|Y_t^{\e, y_0, \sL_{\xi}}|^{p}+\sup_{0\leq t\leq T}\mE|Y_t^{\e, \xi}|^{p}\right\},
\de
and
\ce
\mathscr{J}_4\leq C\e^pT^{p-1}\mE\mathbb{\tilde{E}}\int_0^T\|\g_1(\tilde{X}_s^{\e}, \sL_{X_s^{\e}})\|^{2p}\dif s\leq C\e^pT^p\left\{1+\mE\left(\sup_{0\leq t\leq T}|X_t^{\e}|^{2p}\right)\right\}.
\de

For $\mathscr{J}_5$, using an argument similar to those for $\mathscr{J}_3$ and $\mathscr{J}_4$, we derive
\ce
\mathscr{J}_5\leq C\e^pT^p\left\{1+\mE\left(\sup_{0\leq t\leq T}|X_t^{\e}|^{2p}\right)+\sup_{0\leq t\leq T}\mE|Y_t^{\e, y_0, \sL_{\xi}}|^{2p}+\sup_{0\leq t\leq T}\mE|Y_t^{\e, \xi}|^{2p}\right\}.
\de

For $\mathscr{J}_6$ and $\mathscr{J}_7$, by the Burkholder-Davis-Gundy inequality, the H\"{o}lder inequality and Remark \ref{hgammarem}, it holds that
\ce
&&\mathscr{J}_6+\mathscr{J}_7\no\\
&\leq&C\e^p\mE\left[\int_0^T\|\g_1(X_s^{\e}, \sL_{X_s^{\e}})\|^2\dif s\right]^{\frac p2}+C\e^{\frac p2}\mE\left[\int_0^T\|\g_2(\sL_{X_s^{\e}}, Y_s^{\e, y_0, \sL_{\xi}}, \sL_{Y_s^{\e, \xi}})\|^2\dif s\right]^{\frac p2}\no\\
&\leq&C(\e^p+\e^{\frac p2})T^{\frac p2}\left\{1+\mE\left(\sup_{0\leq t\leq T}|X_t^{\e}|^{p}\right)+\sup_{0\leq t\leq T}\mE|Y_t^{\e, y_0, \sL_{\xi}}|^{p}+\sup_{0\leq t\leq T}\mE|Y_t^{\e, \xi}|^{p}\right\}.
\de

Collecting the above deductions and applying Lemma \ref{xtyt}, we conclude that
\ce
&&\mE\left(\sup_{0\leq t\leq T}\left|\int_0^t\Big(h_1(X_s^{\e}, \sL_{X_s^{\e}},Y_s^{\e, y_0, \sL_{\xi}}, \sL_{Y_s^{\e, \xi}})-\bar{h}_1(X_s^{\e}, \sL_{X_s^{\e}})\Big)\dif s\right|^p\right)\no\\
&\leq& C_T\e^{\frac p2}(1+\mE|\varrho|^{2p}+| y_0|^{2p}+\mE|\xi|^{2p}),
\de
which together with \eqref{xebarxgrow} implies that
$$
\mE\left(\sup_{0 \leq t \leq T}| X_t^{\e}- \bar{X}_t|^p\right)\leq C_T\e^{\frac p2}(1+\mE|\varrho|^{2p}+| y_0|^{2p}+\mE|\xi|^{2p}).
$$
The proof is complete.

\section{Proof of theorem \ref{cltth}}\label{cltthproo}

In this section, we provide the proof of Theorem \ref{cltth}. The first part is devoted to deriving the weak fluctuation estimates. We then complete the proof by combining the result from the previous subsection with the regularity of the solution to the associated Cauchy problem.

\subsection{Weak fluctuation estimates}\label{weakfluct}

In this subsection, we establish the weak fluctuation estimate.

Let $g(t,x,\mu,y,\nu,u,\pi)$ be a function satisfying the centering condition, i.e., for each fixed $(t,x,\mu,u,\pi)\in\mR_+\times\mR^n\times\cP_2(\mR^n)\times\mR^n\times\cP_2(\mR^n)$, the following equation holds:
\be
\int_{\mR^m\times\cP_2(\mR^m)}g(t,x,\mu,y,\nu,u,\pi)(\eta^\mu\times\delta_{\eta^\mu})(\dif y, \dif\nu)=0.
\label{center}
\ee
Then we consider the Poisson equation:
\be
-\cL\phi(t,x,\mu,y,\nu,u,\pi)=g(t,x,\mu,y,\nu,u,\pi),
\label{poissoneq}
\ee
where $(t,x,\mu,u,\pi)$ are treated as parameters. We now state the main result of this subsection.

\bl\label{fluctlem}
Suppose that $(\mathbf{H}_{h_2, \g_2}^2)$ holds with $p\geq2$. Assume that $(\mathbf{H}_{h_1, \g_1}^1)$ and $(\mathbf{H}_{h_2, \g_2}^1)$ hold, $h_2\in C_b^{(1,1),2,(1,1)}(\cP_2(\mR^n)\times\mR^m\times\cP_2(\mR^m),\mR^m), \g_2\in C_b^{(1,1),2,(1,1)}(\cP_2(\mR^n)\times\mR^m\times\cP_2(\mR^m),\mR^{m\times d_2})$.
Let $g\in C_b^{1,2,(1,1),2,(1,1),2,(1,1)}(\mR_+\times\mR^n\times\cP_2(\mR^n)\times\mR^m\times\cP_2(\mR^m)\times\mR^n\times\cP_2(\mR^n),\mR)$ satisfy the centering condition \eqref{center} and assume that $\p_ug(t,\cdot,\mu,y,\nu,u,\pi)\in C_b^1(\mR^n,\mR^n)$. Then there exists a positive constant $C_T$ independent of $\e$ such that
\be
&&\mE\int_0^tg(s,X_s^{\e}, \sL_{X_s^{\e}},Y_s^{\e, y_0, \sL_{\xi}}, \sL_{Y_s^{\e, \xi}},U_s^{\e},\sL_{U_s^{\e}})\dif s\no\\
&\leq&C_T\e+\sqrt\e\mE\int_0^t\left[h_1(X_s^{\e}, \sL_{X_s^{\e}},Y_s^{\e, y_0, \sL_{\xi}}, \sL_{Y_s^{\e, \xi}})-\bar{h}_1(X_s^{\e}, \sL_{X_s^{\e}})\right]\no\\
&&\qquad\qquad\qquad\qquad\cdot\p_u\phi(s,X_s^{\e}, \sL_{X_s^{\e}},Y_s^{\e, y_0, \sL_{\xi}}, \sL_{Y_s^{\e, \xi}},U_s^{\e},\sL_{U_s^{\e}})\dif s\no\\
&&+\sqrt\e\mE\tilde{\mE}\int_0^t\left[h_1(\tilde{X}_s^{\e}, \sL_{X_s^{\e}},\tilde{Y}_s^{\e, y_0, \sL_{\xi}}, \sL_{Y_s^{\e, \xi}})-\bar{h}_1(\tilde{X}_s^{\e}, \sL_{X_s^{\e}})\right]\no\\
&&\qquad\qquad\qquad\qquad\cdot\p_\pi\phi(s,X_s^{\e}, \sL_{X_s^{\e}},Y_s^{\e, y_0, \sL_{\xi}}, \sL_{Y_s^{\e, \xi}},U_s^{\e},\sL_{U_s^{\e}})(\tilde{U}_s^{\e})\dif s,
\label{gfluct}
\ee
where $\phi$ is the solution of the Poisson equation \eqref{poissoneq}.
\el
\begin{proof}
Under the assumptions stated above and according to \cite[Theorem 3.1]{lx}, we have $\phi\in C_b^{1,2,(1,1),2,(1,1),2,(1,1)}(\mR_+\times\mR^n\times\cP_2(\mR^n)\times\mR^m\times\cP_2(\mR^m)\times\mR^n\times\cP_2(\mR^n),\mR)$ and $\p_u\phi(t,\cdot,\mu,y,\nu,u,\pi)\in C_b^1(\mR^n,\mR^n)$.

From the It\^{o} formula, it follows that
\be
&&\phi(t,X_t^{\e}, \sL_{X_t^{\e}},Y_t^{\e, y_0, \sL_{\xi}}, \sL_{Y_t^{\e,\xi}},U_t^{\e},\sL_{U_t^{\e}})-\phi(0,\varrho,\sL_{\varrho},y_0,\sL_{\xi},0,\delta_0)\no\\
&=&\int_0^t\p_s\phi(s,X_s^{\e}, \sL_{X_s^{\e}},Y_s^{\e, y_0, \sL_{\xi}}, \sL_{Y_s^{\e,\xi}},U_s^{\e},\sL_{U_s^{\e}})\dif s\no\\
&&+\int_0^th_1(X_s^{\e}, \sL_{X_s^{\e}},Y_s^{\e, y_0, \sL_{\xi}}, \sL_{Y_s^{\e, \xi}})\cdot\p_x\phi(s,X_s^{\e}, \sL_{X_s^{\e}},Y_s^{\e, y_0, \sL_{\xi}}, \sL_{Y_s^{\e,\xi}},U_s^{\e},\sL_{U_s^{\e}})\dif s\no\\
&&+\int_0^t\p_x\phi(s,X_s^{\e}, \sL_{X_s^{\e}},Y_s^{\e, y_0, \sL_{\xi}}, \sL_{Y_s^{\e,\xi}},U_s^{\e},\sL_{U_s^{\e}})\cdot\g_1(X_s^{\e}, \sL_{X_s^{\e}})\dif B_s\no\\
&&+\frac12\int_0^tTr\left[\g_1\g_1^*(X_s^{\e}, \sL_{X_s^{\e}})\cdot\p_{x}^2\phi(s,X_s^{\e}, \sL_{X_s^{\e}},Y_s^{\e, y_0, \sL_{\xi}}, \sL_{Y_s^{\e,\xi}},U_s^{\e},\sL_{U_s^{\e}})\right]\dif s\no\\
&&+\frac1{\e}\int_0^t h_2(\sL_{X_s^{\e}}, Y_s^{\e, y_0, \sL_{\xi}},\sL_{Y_s^{\e, \xi}})\cdot\p_y\phi(s,X_s^{\e}, \sL_{X_s^{\e}},Y_s^{\e, y_0, \sL_{\xi}}, \sL_{Y_s^{\e,\xi}},U_s^{\e},\sL_{U_s^{\e}})\dif s\no\\
&&+\frac1{\sqrt\e}\int_0^t\p_y\phi(s,X_s^{\e}, \sL_{X_s^{\e}},Y_s^{\e, y_0, \sL_{\xi}}, \sL_{Y_s^{\e,\xi}},U_s^{\e},\sL_{U_s^{\e}})\cdot\g_2(\sL_{X_s^{\e}}, Y_s^{\e, y_0, \sL_{\xi}}, \sL_{Y_s^{\e, \xi}})\dif W_s\no\\
&&+\frac1{2\e}\int_0^tTr\left[\g_2\g_2^*(\sL_{X_s^{\e}}, Y_s^{\e, y_0, \sL_{\xi}}, \sL_{Y_s^{\e, \xi}})\cdot\p_{y}^2\phi(s,X_s^{\e}, \sL_{X_s^{\e}},Y_s^{\e, y_0, \sL_{\xi}}, \sL_{Y_s^{\e,\xi}},U_s^{\e},\sL_{U_s^{\e}})\right]\dif s\no\\
&&+\frac1{\sqrt\e}\int_0^t\left[h_1(X_s^{\e}, \sL_{X_s^{\e}},Y_s^{\e, y_0, \sL_{\xi}}, \sL_{Y_s^{\e, \xi}})-\bar{h}_1(\bar{X}_s, \sL_{\bar{X}_s})\right]\no\\
&&\qquad\qquad\qquad\qquad\cdot\p_u\phi(s,X_s^{\e}, \sL_{X_s^{\e}},Y_s^{\e, y_0, \sL_{\xi}}, \sL_{Y_s^{\e, \xi}},U_s^{\e},\sL_{U_s^{\e}})\dif s\no\\
&&+\frac1{\sqrt\e}\int_0^t\p_u\phi(s,X_s^{\e}, \sL_{X_s^{\e}},Y_s^{\e, y_0, \sL_{\xi}}, \sL_{Y_s^{\e, \xi}},U_s^{\e},\sL_{U_s^{\e}})\cdot[\g_1(X_s^{\e}, \sL_{X_s^{\e}})-\g_1(\bar{X}_s, \sL_{\bar{X}_s})]\dif B_s\no\\
&&+\frac1{2\e}\int_0^tTr\Big[\left(\g_1(X_s^{\e}, \sL_{X_s^{\e}})-\g_1(\bar{X}_s, \sL_{\bar{X}_s})\right)\left(\g_1(X_s^{\e}, \sL_{X_s^{\e}})-\g_1(\bar{X}_s, \sL_{\bar{X}_s})\right)^*\no\\
&&\qquad\qquad\qquad\qquad\cdot\p_{u}^2\phi(s,X_s^{\e}, \sL_{X_s^{\e}},Y_s^{\e, y_0, \sL_{\xi}}, \sL_{Y_s^{\e,\xi}},U_s^{\e},\sL_{U_s^{\e}})\Big]\dif s\no\\
&&+\frac1{\sqrt\e}\int_0^tTr\Big[\g_1(X_s^{\e}, \sL_{X_s^{\e}})\left(\g_1(X_s^{\e}, \sL_{X_s^{\e}})-\g_1(\bar{X}_s, \sL_{\bar{X}_s})\right)^*\no\\
&&\qquad\qquad\qquad\qquad\cdot\p_{x}\p_{u}\phi(s,X_s^{\e}, \sL_{X_s^{\e}},Y_s^{\e, y_0, \sL_{\xi}}, \sL_{Y_s^{\e,\xi}},U_s^{\e},\sL_{U_s^{\e}})\Big]\dif s\no\\
&&+\mathbb{\tilde{E}}\int_0^t h_1(\tilde{X}_s^{\e}, \sL_{X_s^{\e}},\tilde{Y}_s^{\e, y_0, \sL_{\xi}}, \sL_{Y_s^{\e, \xi}})\cdot\p_\mu\phi(s,X_s^{\e}, \sL_{X_s^{\e}},Y_s^{\e, y_0, \sL_{\xi}},
\sL_{Y_s^{\e,\xi}},U_s^{\e},\sL_{U_s^{\e}})(\tilde{X}_s^{\e})\dif s\no\\
&&+\frac12\mathbb{\tilde{E}}\int_0^tTr\left[\g_1\g_1^*(\tilde{X}_s^{\e}, \sL_{X_s^{\e}})\cdot\p_{\tilde{x}}\p_\mu\phi(s,X_s^{\e}, \sL_{X_s^{\e}},Y_s^{\e, y_0, \sL_{\xi}}, \sL_{Y_s^{\e,\xi}},U_s^{\e},\sL_{U_s^{\e}})
(\tilde{X}_s^{\e})\right]\dif s\no\\
&&+\frac1\e\mathbb{\tilde{E}}\int_0^th_2(\sL_{X_s^{\e}}, \tilde{Y}_s^{\e, \tilde{\xi}}, \sL_{Y_s^{\e,\xi}})\cdot\p_\nu\phi(s,X_s^{\e}, \sL_{X_s^{\e}},Y_s^{\e, y_0, \sL_{\xi}}, \sL_{Y_s^{\e,\xi}},U_s^{\e},\sL_{U_s^{\e}})(\tilde{Y}_s^{\e, \tilde{\xi}})\dif s\no\\
&&+\frac1{2\e}\mathbb{\tilde{E}}\int_0^tTr\Big[\g_2\g_2^*(\sL_{X_s^{\e}}, \tilde{Y}_s^{\e, \tilde{\xi}}, \sL_{Y_s^{\e,\xi}})\no\\
&&\qquad\qquad\qquad\qquad\cdot\p_{\tilde{y}}\p_\nu\phi(s,X_s^{\e}, \sL_{X_s^{\e}},Y_s^{\e, y_0, \sL_{\xi}}, \sL_{Y_s^{\e,\xi}},U_s^{\e},\sL_{U_s^{\e}})(\tilde{Y}_s^{\e, \tilde{\xi}})\Big]\dif s\no\\
&&+\frac1{\sqrt\e}\mathbb{\tilde{E}}\int_0^t\left[h_1(\tilde{X}_s^{\e}, \sL_{X_s^{\e}},\tilde{Y}_s^{\e, y_0, \sL_{\xi}}, \sL_{Y_s^{\e, \xi}})-\bar{h}_1(\tilde{\bar{X}}_s, \sL_{\bar{X}_s})\right]\no\\
&&\qquad\qquad\qquad\qquad\cdot\p_\pi\phi(s,X_s^{\e}, \sL_{X_s^{\e}},Y_s^{\e, y_0, \sL_{\xi}}, \sL_{Y_s^{\e, \xi}},U_s^{\e},\sL_{U_s^{\e}})(\tilde{U}_s^{\e})\dif s\no\\
&&+\frac1{2\e}\mathbb{\tilde{E}}\int_0^tTr\Big[\left(\g_1(\tilde{X}_s^{\e}, \sL_{X_s^{\e}})-\g_1(\tilde{\bar{X}}_s, \sL_{\bar{X}_s})\right)\left(\g_1(\tilde{X}_s^{\e}, \sL_{X_s^{\e}})-\g_1(\tilde{\bar{X}}_s,
\sL_{\bar{X}_s})\right)^*\no\\
&&\qquad\qquad\qquad\qquad\cdot\p_{\tilde{u}}\p_\pi\phi(s,X_s^{\e}, \sL_{X_s^{\e}},Y_s^{\e, y_0, \sL_{\xi}}, \sL_{Y_s^{\e,\xi}},U_s^{\e},\sL_{U_s^{\e}})(\tilde{U}_s^{\e})\Big]\dif s,
\label{phiito}
\ee
where the process $(\tilde{X}_t^{\e},\tilde{\bar{X}}_t,\tilde{Y}_t^{\e, y_0, \sL_{\xi}},\tilde{Y}_t^{\e, \tilde{\xi}},\tilde{U}_t^{\e})$ denotes a copy of $(X_t^{\e},\bar{X}_t,Y_t^{\e, y_0, \sL_{\xi}},Y_t^{\e,\xi},U_t^{\e})$ defined on a probability space $(\tilde{\Omega},\tilde{\sF},\{\tilde{\sF}_t\}_{t\geq0},\tilde{\mP})$ which is a copy of $(\Omega,\sF,\{\sF_t\}_{t\geq0},\mP)$, $\tilde{\mE}$ is the expectation taken with respect to $\tilde{\mP}$.
For convenience, we define
\ce
&&\cL_2:=\sum_{i=1}^n[h_{1,i}(x,\mu,y,\nu)-\bar{h}_{1,i}(\bar{x},\bar{\mu})]\frac{\p}{\p_{u_i}},\no\\
&&\cL_3^\e:=\frac1{2\e}\sum_{i,j=1}^n\big([\g_1(x,\mu)-\g_1(\bar{x},\bar{\mu})][\g_1(x,\mu)-\g_1(\bar{x},\bar{\mu})]^*\big)_{ij}\frac{\p^2}{\p_{u_i}\p_{u_j}}\no\\
&&\qquad\quad+\frac1{\sqrt\e}\sum_{i,j=1}^n\big(\g_1(x,\mu)[\g_1(x,\mu)-\g_1(\bar{x},\bar{\mu})]^*\big)_{ij}\frac{\p^2}{\p_{x_i}\p_{u_j}}.
\de
Then Eq.\eqref{phiito} can be rewritten as
\be
&&\phi(t,X_t^{\e}, \sL_{X_t^{\e}},Y_t^{\e, y_0, \sL_{\xi}}, \sL_{Y_t^{\e,\xi}},U_t^{\e},\sL_{U_t^{\e}})-\phi(0,\varrho,\sL_{\varrho},y_0,\sL_{\xi},0,\delta_0)\no\\
&=&\frac1\e\int_0^t\cL\phi(s,X_s^{\e}, \sL_{X_s^{\e}},Y_s^{\e, y_0, \sL_{\xi}}, \sL_{Y_s^{\e,\xi}},U_s^{\e},\sL_{U_s^{\e}})\dif s\no\\
&&+\int_0^t(\p_s+{\cL}_1+{\cL}_3^\e)\phi(s,X_s^{\e}, \sL_{X_s^{\e}},Y_s^{\e, y_0, \sL_{\xi}}, \sL_{Y_s^{\e,\xi}},U_s^{\e},\sL_{U_s^{\e}})\dif s\no\\
&&+\frac1{\sqrt\e}\int_0^t{\cL}_2\phi(s,X_s^{\e}, \sL_{X_s^{\e}},Y_s^{\e, y_0, \sL_{\xi}}, \sL_{Y_s^{\e,\xi}},U_s^{\e},\sL_{U_s^{\e}})\dif s\no\\
&&+\mathbb{\tilde{E}}\int_0^t h_1(\tilde{X}_s^{\e}, \sL_{X_s^{\e}},\tilde{Y}_s^{\e, y_0, \sL_{\xi}}, \sL_{Y_s^{\e, \xi}})\cdot\p_\mu\phi(s,X_s^{\e}, \sL_{X_s^{\e}},Y_s^{\e, y_0, \sL_{\xi}}, \sL_{Y_s^{\e,\xi}},
U_s^{\e},\sL_{U_s^{\e}})(\tilde{X}_s^{\e})\dif s\no\\
&&+\frac12\mathbb{\tilde{E}}\int_0^tTr\left[\g_1\g_1^*(\tilde{X}_s^{\e}, \sL_{X_s^{\e}})\cdot\p_{\tilde{x}}\p_\mu\phi(s,X_s^{\e}, \sL_{X_s^{\e}},Y_s^{\e, y_0, \sL_{\xi}}, \sL_{Y_s^{\e,\xi}},U_s^{\e},\sL_{U_s^{\e}})
(\tilde{X}_s^{\e})\right]\dif s\no\\
&&+\frac1{\sqrt\e}\mathbb{\tilde{E}}\int_0^t\left[h_1(\tilde{X}_s^{\e}, \sL_{X_s^{\e}},\tilde{Y}_s^{\e, y_0, \sL_{\xi}}, \sL_{Y_s^{\e, \xi}})-\bar{h}_1(\tilde{\bar{X}}_s, \sL_{\bar{X}_s})\right]\no\\
&&\qquad\qquad\qquad\qquad\cdot\p_\pi\phi(s,X_s^{\e}, \sL_{X_s^{\e}},Y_s^{\e, y_0, \sL_{\xi}}, \sL_{Y_s^{\e, \xi}},U_s^{\e},\sL_{U_s^{\e}})(\tilde{U}_s^{\e})\dif s\no\\
&&+\frac1{2\e}\mathbb{\tilde{E}}\int_0^tTr\Big[\left(\g_1(\tilde{X}_s^{\e}, \sL_{X_s^{\e}})-\g_1(\tilde{\bar{X}}_s, \sL_{\bar{X}_s})\right)\left(\g_1(\tilde{X}_s^{\e}, \sL_{X_s^{\e}})-\g_1(\tilde{\bar{X}}_s, \sL_{\bar{X}_s})\right)^*\no\\
&&\qquad\qquad\qquad\qquad\cdot\p_{\tilde{u}}\p_\pi\phi(s,X_s^{\e}, \sL_{X_s^{\e}},Y_s^{\e, y_0, \sL_{\xi}}, \sL_{Y_s^{\e,\xi}},U_s^{\e},\sL_{U_s^{\e}})(\tilde{U}_s^{\e})\Big]\dif s\no\\
&&+\int_0^t\p_x\phi(s,X_s^{\e}, \sL_{X_s^{\e}},Y_s^{\e, y_0, \sL_{\xi}}, \sL_{Y_s^{\e,\xi}},U_s^{\e},\sL_{U_s^{\e}})\cdot\g_1(X_s^{\e}, \sL_{X_s^{\e}})\dif B_s\no\\
&&+\frac1{\sqrt\e}\int_0^t\p_y\phi(s,X_s^{\e}, \sL_{X_s^{\e}},Y_s^{\e, y_0, \sL_{\xi}}, \sL_{Y_s^{\e,\xi}},U_s^{\e},\sL_{U_s^{\e}})\cdot\g_2(\sL_{X_s^{\e}}, Y_s^{\e, y_0, \sL_{\xi}}, \sL_{Y_s^{\e, \xi}})\dif W_s\no\\
&&+\frac1{\sqrt\e}\int_0^t\p_u\phi(s,X_s^{\e}, \sL_{X_s^{\e}},Y_s^{\e, y_0, \sL_{\xi}}, \sL_{Y_s^{\e, \xi}},U_s^{\e},\sL_{U_s^{\e}})\cdot[\g_1(X_s^{\e}, \sL_{X_s^{\e}})-\g_1(\bar{X}_s, \sL_{\bar{X}_s})]\dif B_s.\no\\
\label{phiitorew}
\ee
By multiplying both sides of \eqref{phiitorew} by $\e$ and taking expectations, and then applying \eqref{poissoneq}, we obtain
\ce
&&\mE\int_0^tg(s,X_s^{\e}, \sL_{X_s^{\e}},Y_s^{\e, y_0, \sL_{\xi}}, \sL_{Y_s^{\e, \xi}},U_s^{\e},\sL_{U_s^{\e}})\dif s\no\\
&=&-\mE\int_0^t\cL\phi(s,X_s^{\e}, \sL_{X_s^{\e}},Y_s^{\e, y_0, \sL_{\xi}}, \sL_{Y_s^{\e,\xi}},U_s^{\e},\sL_{U_s^{\e}})\dif s\no\\
&=&\e\mE\left[\phi(0,\varrho,\sL_{\varrho},y_0,\sL_{\xi},0,\delta_0)-\phi(t,X_t^{\e}, \sL_{X_t^{\e}},Y_t^{\e, y_0, \sL_{\xi}}, \sL_{Y_t^{\e,\xi}},U_t^{\e},\sL_{U_t^{\e}})\right]\no\\
&&+\e\mE\int_0^t(\p_s+{\cL}_1+{\cL}_3^\e)\phi(s,X_s^{\e}, \sL_{X_s^{\e}},Y_s^{\e, y_0, \sL_{\xi}}, \sL_{Y_s^{\e,\xi}},U_s^{\e},\sL_{U_s^{\e}})\dif s\no\\
&&+\sqrt\e\mE\int_0^t{\cL}_2\phi(s,X_s^{\e}, \sL_{X_s^{\e}},Y_s^{\e, y_0, \sL_{\xi}}, \sL_{Y_s^{\e,\xi}},U_s^{\e},\sL_{U_s^{\e}})\dif s\no\\
&&+\e\mE\mathbb{\tilde{E}}\int_0^t h_1(\tilde{X}_s^{\e}, \sL_{X_s^{\e}},\tilde{Y}_s^{\e, y_0, \sL_{\xi}}, \sL_{Y_s^{\e, \xi}})\no\\
&&\qquad\qquad\qquad\qquad\cdot\p_\mu\phi(s,X_s^{\e}, \sL_{X_s^{\e}},Y_s^{\e, y_0, \sL_{\xi}}, \sL_{Y_s^{\e,\xi}},U_s^{\e},\sL_{U_s^{\e}})(\tilde{X}_s^{\e})\dif s\no\\
&&+\frac\e2\mE\mathbb{\tilde{E}}\int_0^tTr\left[\g_1\g_1^*(\tilde{X}_s^{\e}, \sL_{X_s^{\e}})\cdot\p_{\tilde{x}}\p_\mu\phi(s,X_s^{\e}, \sL_{X_s^{\e}},Y_s^{\e, y_0, \sL_{\xi}}, \sL_{Y_s^{\e,\xi}},U_s^{\e},\sL_{U_s^{\e}})
(\tilde{X}_s^{\e})\right]\dif s\no\\
&&+\sqrt\e\mE\mathbb{\tilde{E}}\int_0^t\left[h_1(\tilde{X}_s^{\e}, \sL_{X_s^{\e}},\tilde{Y}_s^{\e, y_0, \sL_{\xi}}, \sL_{Y_s^{\e, \xi}})-\bar{h}_1(\tilde{\bar{X}}_s, \sL_{\bar{X}_s})\right]\no\\
&&\qquad\qquad\qquad\qquad\cdot\p_\pi\phi(s,X_s^{\e}, \sL_{X_s^{\e}},Y_s^{\e, y_0, \sL_{\xi}}, \sL_{Y_s^{\e, \xi}},U_s^{\e},\sL_{U_s^{\e}})(\tilde{U}_s^{\e})\dif s\no\\
&&+\frac12\mE\mathbb{\tilde{E}}\int_0^tTr\Big[\left(\g_1(\tilde{X}_s^{\e}, \sL_{X_s^{\e}})-\g_1(\tilde{\bar{X}}_s, \sL_{\bar{X}_s})\right)\left(\g_1(\tilde{X}_s^{\e}, \sL_{X_s^{\e}})-\g_1(\tilde{\bar{X}}_s, \sL_{\bar{X}_s})\right)^*\no\\
&&\qquad\qquad\qquad\qquad\cdot\p_{\tilde{u}}\p_\pi\phi(s,X_s^{\e}, \sL_{X_s^{\e}},Y_s^{\e, y_0, \sL_{\xi}}, \sL_{Y_s^{\e,\xi}},U_s^{\e},\sL_{U_s^{\e}})(\tilde{U}_s^{\e})\Big]\dif s\no\\
&=&\sum_{i=1}^7I_i^\e.
\de

For $I_1^\e$ and $I_2^\e$, the regularity of $\phi$ together with Lemma \ref{xtyt} and Theorem \ref{xbarxp} yields that
\ce
I_1^\e\leq C\e\big[T+(\mE|X_t^{\e}-\varrho|^2)^{\frac12}+\mE|Y_t^{\e, y_0, \sL_{\xi}}-y_0|+(\mE|Y_t^{\e,\xi}-\xi|^2)^{\frac12}+(\mE|U_t^{\e}|^2)^{\frac12}\big]\leq C_T\e,
\de
and
\ce
I_2^\e&\leq&C\e\mE\int_0^t\Big[1+|h_1(X_s^{\e}, \sL_{X_s^{\e}},Y_s^{\e, y_0, \sL_{\xi}}, \sL_{Y_s^{\e, \xi}})|+\|\g_1(X_s^{\e}, \sL_{X_s^{\e}})\|^2\Big]\dif s\\
&&+C\e\mE\int_0^t\Big[\frac1\e\|\g_1(X_s^{\e}, \sL_{X_s^{\e}})-\g_1(\bar{X}_s, \sL_{\bar{X}_s})\|^2\\
&&\qquad\qquad\qquad+\frac1{\sqrt\e}\|\g_1(X_s^{\e}, \sL_{X_s^{\e}})\|\cdot\|\g_1(X_s^{\e}, \sL_{X_s^{\e}})-\g_1(\bar{X}_s, \sL_{\bar{X}_s})\|\Big]\dif s\\
&\leq&C\e\int_0^t[1+\mE|X_s^{\e}|^2+\mE|Y_s^{\e, y_0, \sL_{\xi}}|^2+\mE|Y_s^{\e, \xi}|^2+\mE|U_s^{\e}|^2]\dif s\\
&\leq&C_T\e.
\de

Similar arguments to that for $I_1^\e$ and $I_2^\e$ imply that
\ce
I_4^\e+I_5^\e+I_7^\e\leq C_T\e.
\de

For $I_3^\e$, by relying on the regularity of $\phi$, along with \eqref{barh1lip} and Theorem \ref{xbarxp}, one can deduce that
\ce
I_3^\e
&=&\sqrt\e\mE\int_0^t\Big[h_1(X_s^{\e}, \sL_{X_s^{\e}},Y_s^{\e, y_0, \sL_{\xi}}, \sL_{Y_s^{\e, \xi}})-\bar{h}_1(X_s^{\e}, \sL_{X_s^{\e}})+\bar{h}_1(X_s^{\e}, \sL_{X_s^{\e}})-\bar{h}_1(\bar{X}_s, \sL_{\bar{X}_s})\Big]\\
&&\qquad\qquad\qquad\qquad\cdot\p_u\phi(s,X_s^{\e}, \sL_{X_s^{\e}},Y_s^{\e, y_0, \sL_{\xi}}, \sL_{Y_s^{\e, \xi}},U_s^{\e},\sL_{U_s^{\e}})\dif s\\
&\leq&C_T\e+\sqrt\e\mE\int_0^t\Big[h_1(X_s^{\e}, \sL_{X_s^{\e}},Y_s^{\e, y_0, \sL_{\xi}}, \sL_{Y_s^{\e, \xi}})-\bar{h}_1(X_s^{\e}, \sL_{X_s^{\e}})\Big]\\
&&\qquad\qquad\qquad\qquad\qquad\qquad\cdot\p_u\phi(s,X_s^{\e}, \sL_{X_s^{\e}},Y_s^{\e, y_0, \sL_{\xi}}, \sL_{Y_s^{\e, \xi}},U_s^{\e},\sL_{U_s^{\e}})\dif s.\\
\de

Similarly, for $I_6^\e$, it follows that
\ce
I_6^\e
&\leq&C_T\e+\sqrt\e\mE\mathbb{\tilde{E}}\int_0^t\left[h_1(\tilde{X}_s^{\e}, \sL_{X_s^{\e}},\tilde{Y}_s^{\e, y_0, \sL_{\xi}}, \sL_{Y_s^{\e, \xi}})-\bar{h}_1(\tilde{X}_s^{\e}, \sL_{X_s^{\e}})\right]\\
&&\qquad\qquad\qquad\qquad\qquad\cdot\p_\pi\phi(s,X_s^{\e}, \sL_{X_s^{\e}},Y_s^{\e, y_0, \sL_{\xi}}, \sL_{Y_s^{\e, \xi}},U_s^{\e},\sL_{U_s^{\e}})(\tilde{U}_s^{\e})\dif s.
\de

Combining the above computations, we conclude that the inequality \eqref{gfluct} holds. The proof is complete.
\end{proof}

\br
We mention that by Lemma \ref{xtyt} and the regularity of $\phi$, it holds that
\be
\left|\mE\int_0^tg(s,X_s^{\e}, \sL_{X_s^{\e}},Y_s^{\e, y_0, \sL_{\xi}}, \sL_{Y_s^{\e,\xi}},U_s^{\e},\sL_{U_s^{\e}})\dif s\right|\leq C_T\sqrt\e.
\label{gfluctrem}
\ee
However, since the terms involving expectations in \eqref{gfluct} will be important for investigating the central limit theorem, we retain them for future use.
\label{fluctrem}
\er

\subsection{Proof of Theorem \ref{cltth}}\label{cltthproode}
In this subsection, we prove Theorem \ref{cltth}.

Let $\bar{X}^{s,\varrho}_t$ be the unique solution to Eq.\eqref{barxeq} with initial data $\varrho\in L^2(\Omega, \sF_0, \mP; \mR^n)$ at time $s$, and let $U^{s,\varrho,\kappa}_t$ be the unique solution to Eq.\eqref{uteq} with the initial value $\kappa\in L^2(\Omega, \sF_0, \mP; \mR^n)$ at time $s$. Specifically, for $t\geq s$,
\ce
\left\{\begin{array}{l}
\dif\bar{X}^{s,\varrho}_t=\bar{h}_1(\bar{X}^{s,\varrho}_t, \sL_{\bar{X}^{s,\varrho}_t})\dif t+\g_1(\bar{X}^{s,\varrho}_t, \sL_{\bar{X}^{s,\varrho}_t})\dif B_t, \\
\bar{X}^{s,\varrho}_s=\varrho,
\end{array}
\right.
\de
and
\ce
\left\{\begin{array}{l}
\dif U^{s,\varrho,\kappa}_t=\p_x\bar{h}_1(\bar{X}^{s,\varrho}_t, \sL_{\bar{X}^{s,\varrho}_t})U^{s,\varrho,\kappa}_t\dif t+\tilde{\mE}[\p_\mu\bar{h}_1(\bar{X}^{s,\varrho}_t, \sL_{\bar{X}^{s,\varrho}_t})
(\tilde{\bar{X}}^{s,\tilde{\varrho}}_t)\tilde{U}^{s,\tilde{\varrho},\tilde{\kappa}}_t]\dif t\\
\qquad\qquad+\p_x\g_1(\bar{X}^{s,\varrho}_t, \sL_{\bar{X}^{s,\varrho}_t})U^{s,\varrho,\kappa}_t\dif B_t+\tilde{\mE}[\p_\mu\g_1(\bar{X}^{s,\varrho}_t, \sL_{\bar{X}^{s,\varrho}_t})(\tilde{\bar{X}}^{s,\tilde{\varrho}}_t)
\tilde{U}^{s,\tilde{\varrho},\tilde{\kappa}}_t]\dif B_t\\
\qquad\qquad+\Upsilon(\bar{X}^{s,\varrho}_t, \sL_{\bar{X}^{s,\varrho}_t})\dif V_t,\\
U^{s,\varrho,\kappa}_s=\kappa.
\end{array}
\right.
\de
Given a fixed $T>0$ and a function $\psi: \cP_2(\mR^n)\rightarrow\mR$, we consider the following Cauchy problem on $[0,T]\times\cP_2(\mR^n)\times\cP_2(\mR^n)$ :
\be
\left\{\begin{array}{l}
\p_tf(t,\sL_{\varrho},\sL_{\kappa})+\mE[\bar{h}_1(\varrho,\sL_{\varrho})\cdot\p_\mu f(t,\sL_{\varrho},\sL_{\kappa})(\varrho)]\\
\quad+\frac12\mE\left[Tr\big(\g_1\g_1^*(\varrho,\sL_{\varrho})\cdot\p_x\p_\mu f(t,\sL_{\varrho},\sL_{\kappa})(\varrho)\big)\right]\\
\quad+\mE\left[\left(\p_x\bar{h}_1(\varrho,\sL_{\varrho})\kappa+\tilde{\mE}[\p_\mu\bar{h}_1(\varrho,\sL_{\varrho})(\tilde{\varrho})\tilde{\kappa}]\right)\cdot\p_\pi f(t,\sL_{\varrho},\sL_{\kappa})(\kappa)\right]\\
\quad+\frac12\mE\Big[Tr\Big(\big[\p_x\g_1(\varrho,\sL_{\varrho})\kappa+\tilde{\mE}[\p_\mu\g_1(\varrho,\sL_{\varrho})(\tilde{\varrho})\tilde{\kappa}]\big]\big[\p_x\g_1(\varrho,\sL_{\varrho})\kappa
+\tilde{\mE}[\p_\mu\g_1(\varrho,\sL_{\varrho})(\tilde{\varrho})\tilde{\kappa}]\big]^*\\
\qquad\qquad\qquad\qquad\cdot\p_u\p_\pi f(t,\sL_{\varrho},\sL_{\kappa})(\kappa)\Big)\Big]\\
\quad+\mE\Big[Tr\Big(\overline{\chi_{h_1}\Psi^*}(\varrho,\sL_{\varrho})\cdot\p_u\p_\pi f(t,\sL_{\varrho},\sL_{\kappa})(\kappa)\Big)\Big]=0,\\
f(T,\sL_{\varrho},\sL_{\kappa})=\psi(\sL_{\kappa}),
\end{array}
\right.
\label{cauchypro}
\ee
where $\chi_{h_1}(x,\mu,y,\nu):=h_1(x,\mu,y,\nu)-\bar{h}_1(x,\mu)$. Then under the assumptions in Theorem \ref{cltth}, according to \cite[Corollary 2.4]{lwx}, we have $\bar{h}_1\in\big(C_b^{4,(1,3)}\cap C_b^{4,(2,2)}\cap C_b^{4,(3,1)}\big)(\mR^n\times\cP_2(\mR^n),\mR^n)$. Recall that $\Psi$ solves \eqref{h1poieq}. By \cite[Theorem 2.3]{lwx}, we obtain $\Psi\in \big(C_b^{4,(3,1),6,(3,3)}\cap \mathbb{C}_b^{4,(1,3),4,(2,2)}\cap\mathbb{C}_b^{4,(2,2),2,(1,1)}\big)(\mR^n\times\cP_2(\mR^n)
\times\mR^m\times\cP_2(\mR^m),\mR^n)$. Then $\overline{\chi_{h_1}\Psi^*}\in\big(C_b^{4,(1,3)}\cap C_b^{4,(2,2)}\cap C_b^{4,(3,1)}\big)(\mR^n\times\cP_2(\mR^n),\mR)$. Following a similar argument as in \cite[Theorem 7.2]{blpr}, it can be concluded that there exists a unique solution $f\in C_b^{1,(2,1),(3,1)}([0,T]\times\cP_2(\mR^n)\times\cP_2(\mR^n),\mR)$ to Eq.\eqref{cauchypro}, which is given by
\be
f(t,\sL_{\varrho},\sL_{\kappa}):=\psi(\sL_{U^{t,\varrho,\kappa}_T}).
\label{cauchysol}
\ee
Moreover, we deduce that the mapping $(t,\mu,u,\pi)\rightarrow\p_\pi f(t,\mu,\pi)(u)$ is in $C_b^{1,(1,1),2,(1,1)}([0,T]\times\cP_2(\mR^n)\times\mR^n\times\cP_2(\mR^n),\mR^n)$; the mapping $(t,\mu,u,\pi)\rightarrow\p_u\p_\pi f(t,\mu,\pi)(u)$ is in $C_b^{1,(1,1),2,(1,1)}([0,T]\times\cP_2(\mR^n)\times\mR^n\times\cP_2(\mR^n),\mR^{n\times n})$; the mapping $(t,\mu,u,\tilde{u},\pi)\rightarrow\p_\pi^2 f(t,\mu,\pi)(u,\tilde{u})$ is in $C_b^{1,(1,1),2,2,(1,1)}$
$([0,T]\times\cP_2(\mR^n)\times\mR^n\times\mR^n\times\cP_2(\mR^n),\mR^{n\times n})$; and the mapping $(t,x,\mu,\pi)\rightarrow\p_\mu f(t,\mu,\pi)(x)$ is in $C_b^{1,2,(1,1),(1,1)}([0,T]\times\mR^n\times\cP_2(\mR^n)\times\cP_2(\mR^n),\mR^n)$.

Besides, employing a deduction similar to that in \cite[Theorem 3.2]{xq}, we deduce that the solution $U^{\e}$ of Eq.\eqref{ueeq} weakly converges to the solution $U$ of Eq.\eqref{uteq} in $C\left([0, T] ; \mR^{n}\right)$, as $\e\rightarrow0$.

Next, we focus on proving the estimate \eqref{weakrate}. Let $f(t,\sL_{\varrho},\sL_{\kappa})$ be defined by \eqref{cauchysol}. Then we deduce that
\ce
\Theta(\e):&=&\psi(\sL_{U^\e_T})-\psi(\sL_{U_T})=\psi(\sL_{U^{T,X_T^\e,U_T^\e}_T})-\psi(\sL_{U^{0,\varrho,0}_T})\\
&=&f(T,\sL_{X_T^\e},\sL_{U_T^\e})-f(0,\sL_{\varrho},\delta_0).
\de

Applying the It\^{o} formula to $f(t,\sL_{X_t^\e},\sL_{U^\e_t})$, we derive that
\be
\Theta(\e)&=&\mE\int_0^T\Big[\p_tf(t,\sL_{X_t^\e},\sL_{U^\e_t})+h_1(X_t^{\e}, \sL_{X_t^{\e}},Y_t^{\e, y_0, \sL_{\xi}}, \sL_{Y_t^{\e, \xi}})\cdot\p_\mu f(t,\sL_{X_t^\e},\sL_{U^\e_t})(X_t^\e)\no\\
&&\qquad\qquad+\frac12Tr\Big(\g_1\g_1^*(X_t^{\e}, \sL_{X_t^{\e}})\cdot\p_x\p_\mu f(t,\sL_{X_t^\e},\sL_{U^\e_t})(X_t^\e)\Big)\Big]\dif t\no\\
&&+\frac1{\sqrt\e}\mE\int_0^T\left[\chi_{h_1}(X_t^{\e}, \sL_{X_t^{\e}},Y_t^{\e, y_0, \sL_{\xi}}, \sL_{Y_t^{\e, \xi}})+\bar{h}_1(X_t^{\e}, \sL_{X_t^{\e}})-\bar{h}_1(\bar{X}_t, \sL_{\bar{X}_t})\right]\no\\
&&\qquad\qquad\qquad\qquad\qquad\qquad\qquad\qquad\cdot\p_\pi f(t,\sL_{X_t^\e},\sL_{U^\e_t})(U_t^\e)\dif t\no\\
&&+\frac1{2\e}\mE\int_0^TTr\Big(\left[\g_1(X_t^{\e}, \sL_{X_t^{\e}})-\g_1(\bar{X}_t, \sL_{\bar{X}_t})\right]\left[\g_1(X_t^{\e}, \sL_{X_t^{\e}})-\g_1(\bar{X}_t, \sL_{\bar{X}_t})\right]^*\no\\
&&\qquad\qquad\qquad\qquad\qquad\qquad\qquad\cdot\p_u\p_\pi f(t,\sL_{X_t^\e},\sL_{U^\e_t})(U_t^\e)\Big)\dif t.
\label{fito}
\ee

By the definition of $\bar{h}_1(x,\mu)$, the function $\chi_{h_1}(x,\mu,y,\nu)\cdot\p_\pi f(t,\mu,\pi)(u)$ satisfies the centering condition given in \eqref{center}. Recall that $\Psi(x,\mu,y,\nu)$ is the solution to the Poisson equation \eqref{h1poieq}. Then we introduce the function
\ce
\hat{\Psi}(t,x,\mu,y,\nu,u,\pi)=\Psi(x,\mu,y,\nu)\cdot\p_\pi f(t,\mu,\pi)(u),
\de
which yields that
\ce
-\cL\hat{\Psi}(t,x,\mu,y,\nu,u,\pi)=\chi_{h_1}(x,\mu,y,\nu)\cdot\p_\pi f(t,\mu,\pi)(u).
\de
Moreover, we have $\chi_{h_1}\cdot\p_\pi f\in C_b^{1,2,(1,1),2,(1,1),2,(1,1)}(\mR_+\times\mR^n\times\cP_2(\mR^n)\times\mR^m\times\cP_2(\mR^m)\times\mR^n\times\cP_2(\mR^n),\mR)$. Thus, combining these with \eqref{gfluct}, we obtain that
\be
&&\frac1{\sqrt\e}\mE\int_0^T\chi_{h_1}(X_t^{\e}, \sL_{X_t^{\e}},Y_t^{\e, y_0, \sL_{\xi}}, \sL_{Y_t^{\e, \xi}})\cdot\p_\pi f(t,\sL_{X_t^\e},\sL_{U^\e_t})(U_t^\e)\dif t\no\\
&\leq&C_T\sqrt\e+\mE\int_0^T\chi_{h_1}(X_t^{\e}, \sL_{X_t^{\e}},Y_t^{\e, y_0, \sL_{\xi}}, \sL_{Y_t^{\e, \xi}})\no\\
&&\qquad\qquad\qquad\qquad\cdot\p_u\hat{\Psi}(t,X_t^{\e}, \sL_{X_t^{\e}},Y_t^{\e, y_0, \sL_{\xi}}, \sL_{Y_t^{\e, \xi}},U_t^{\e},\sL_{U_t^{\e}})\dif t\no\\
&&+\mE\tilde{\mE}\int_0^T\chi_{h_1}(\tilde{X}_t^{\e}, \sL_{X_t^{\e}},\tilde{Y}_t^{\e, y_0, \sL_{\xi}}, \sL_{Y_t^{\e, \xi}})\no\\
&&\qquad\qquad\qquad\qquad\cdot\p_\pi\hat{\Psi}(t,X_t^{\e}, \sL_{X_t^{\e}},Y_t^{\e, y_0, \sL_{\xi}}, \sL_{Y_t^{\e, \xi}},U_t^{\e},\sL_{U_t^{\e}})(\tilde{U}_t^{\e})\dif t.
\label{chih1papif}
\ee
For the last term on the right hand side of the above inequality, it is straightforward to observe that for any fixed $(x,y,u)\in\mR^n\times\mR^m\times\mR^n$, the function
\ce
(t,\tilde{x},\mu,\tilde{y},\nu,\tilde{u},\pi)\rightarrow\chi_{h_1}(\tilde{x},\mu,\tilde{y},\nu)\cdot\p_\pi\hat{\Psi}(t,x,\mu,y,\nu,u,\pi)(\tilde{u})
\de
satisfies the centering condition \eqref{center}. Hence, by \eqref{gfluctrem}, it follows that
\ce
&&\mE\tilde{\mE}\int_0^T\chi_{h_1}(\tilde{X}_t^{\e}, \sL_{X_t^{\e}},\tilde{Y}_t^{\e, y_0, \sL_{\xi}}, \sL_{Y_t^{\e, \xi}})\cdot\p_\pi\hat{\Psi}(t,X_t^{\e}, \sL_{X_t^{\e}},Y_t^{\e, y_0, \sL_{\xi}}, \sL_{Y_t^{\e, \xi}},U_t^{\e},\sL_{U_t^{\e}})(\tilde{U}_t^{\e})\dif t\no\\
&\leq&C_T\sqrt\e,
\de
which together with \eqref{fito} and \eqref{chih1papif} implies that
\ce
\Theta(\e)
&\leq&C_T\sqrt\e+\int_0^T\p_tf(t,\sL_{X_t^\e},\sL_{U^\e_t})\dif t\no\\
&&+\mE\int_0^Th_1(X_t^{\e}, \sL_{X_t^{\e}},Y_t^{\e, y_0, \sL_{\xi}}, \sL_{Y_t^{\e, \xi}})\cdot\p_\mu f(t,\sL_{X_t^\e},\sL_{U^\e_t})(X_t^\e)\dif t\no\\
&&+\frac12\mE\int_0^TTr\Big(\g_1\g_1^*(X_t^{\e}, \sL_{X_t^{\e}})\cdot\p_x\p_\mu f(t,\sL_{X_t^\e},\sL_{U^\e_t})(X_t^\e)\Big)\dif t\no\\
&&+\frac1{\sqrt\e}\mE\int_0^T\left[\bar{h}_1(X_t^{\e}, \sL_{X_t^{\e}})-\bar{h}_1(\bar{X}_t, \sL_{\bar{X}_t})\right]\cdot\p_\pi f(t,\sL_{X_t^\e},\sL_{U^\e_t})(U_t^\e)\dif t\no\\
&&+\mE\int_0^T\chi_{h_1}(X_t^{\e}, \sL_{X_t^{\e}},Y_t^{\e, y_0, \sL_{\xi}}, \sL_{Y_t^{\e, \xi}})\no\\
&&\qquad\qquad\qquad\qquad\cdot\p_u\hat{\Psi}(t,X_t^{\e}, \sL_{X_t^{\e}},Y_t^{\e, y_0, \sL_{\xi}}, \sL_{Y_t^{\e, \xi}},U_t^{\e},\sL_{U_t^{\e}})\dif t\no\\
&&+\frac1{2\e}\mE\int_0^TTr\Big(\left[\g_1(X_t^{\e}, \sL_{X_t^{\e}})-\g_1(\bar{X}_t, \sL_{\bar{X}_t})\right]\left[\g_1(X_t^{\e}, \sL_{X_t^{\e}})-\g_1(\bar{X}_t, \sL_{\bar{X}_t})\right]^*\no\\
&&\qquad\qquad\qquad\qquad\qquad\qquad\qquad\cdot\p_u\p_\pi f(t,\sL_{X_t^\e},\sL_{U^\e_t})(U_t^\e)\Big)\dif t.
\de
By integrating this with \eqref{cauchypro}, one can conclude that
\ce
\Theta(\e)
&\leq&C_T\sqrt\e+\mE\int_0^T\chi_{h_1}(X_t^{\e}, \sL_{X_t^{\e}},Y_t^{\e, y_0, \sL_{\xi}}, \sL_{Y_t^{\e, \xi}})\cdot\p_\mu f(t,\sL_{X_t^\e},\sL_{U^\e_t})(X_t^\e)\dif t\no\\
&&+\frac1{\sqrt\e}\mE\int_0^T\Big[\bar{h}_1(X_t^{\e}, \sL_{X_t^{\e}})-\bar{h}_1(\bar{X}_t, \sL_{\bar{X}_t})-\sqrt\e\p_x\bar{h}_1(X_t^{\e}, \sL_{X_t^{\e}})U^\e_t\\
&&\qquad\qquad\qquad-\sqrt\e\tilde{\mE}[\p_\mu\bar{h}_1(X_t^{\e}, \sL_{X_t^{\e}})(\tilde{X}_t^{\e})\tilde{U}_t^{\e}]\Big]\cdot\p_\pi f(t,\sL_{X_t^\e},\sL_{U^\e_t})(U_t^\e)\dif t\no\\
&&+\frac1{2\e}\mE\int_0^TTr\Big(\Big[\left[\g_1(X_t^{\e}, \sL_{X_t^{\e}})-\g_1(\bar{X}_t, \sL_{\bar{X}_t})\right]\left[\g_1(X_t^{\e}, \sL_{X_t^{\e}})-\g_1(\bar{X}_t, \sL_{\bar{X}_t})\right]^*\no\\
&&\qquad\qquad\qquad\quad-\big[\sqrt\e\p_x\g_1(X_t^{\e}, \sL_{X_t^{\e}})U^\e_t+\sqrt\e\tilde{\mE}[\p_\mu\g_1(X_t^{\e}, \sL_{X_t^{\e}})(\tilde{X}_t^{\e})\tilde{U}_t^{\e}]\big]\\
&&\qquad\qquad\qquad\qquad\big[\sqrt\e\p_x\g_1(X_t^{\e}, \sL_{X_t^{\e}})U^\e_t+\sqrt\e\tilde{\mE}[\p_\mu\g_1(X_t^{\e}, \sL_{X_t^{\e}})(\tilde{X}_t^{\e})\tilde{U}_t^{\e}]\big]^*\Big]\\
&&\qquad\qquad\qquad\qquad\qquad\cdot\p_u\p_\pi f(t,\sL_{X_t^\e},\sL_{U^\e_t})(U_t^\e)\Big)\dif t\\
&&+\mE\int_0^TTr\Big(\Big[\chi_{h_1}\Psi^*(X_t^{\e}, \sL_{X_t^{\e}},Y_t^{\e, y_0, \sL_{\xi}}, \sL_{Y_t^{\e, \xi}})-\overline{\chi_{h_1}\Psi^*}(X_t^{\e}, \sL_{X_t^{\e}})\Big]\\
&&\qquad\qquad\qquad\qquad\qquad\cdot\p_u\p_\pi f(t,\sL_{X_t^\e},\sL_{U^\e_t})(U_t^\e)\Big)\dif t\\
&=&C_T\sqrt\e+\sum_{i=1}^4J_i^\e.
\de

For $J_1^\e$, observing that the function
\ce
(t,x,\mu,y,\nu,\pi)\rightarrow\chi_{h_1}(x,\mu,y,\nu)\cdot\p_\mu f(t,\mu,\pi)(x)
\de
satisfies the centering condition \eqref{center} and belongs to $C_b^{1,2,(1,1),2,(1,1),(1,1)}(\mR_+\times\mR^n\times\cP_2(\mR^n)\times\mR^m\times\cP_2(\mR^m)\times\cP_2(\mR^n),\mR)$. By \eqref{gfluctrem}, it holds that
\ce
J_1^\e\leq C_T\sqrt\e.
\de

Similarly to $J_1^\e$, we arrive at
\ce
J_4^\e\leq C_T\sqrt\e.
\de

For $J_2^\e$, the mean value theorem (cf. \cite[Lemma 4.1]{kn}) together with Theorem \ref{xbarxp}, implies that
\ce
J_2^\e&\leq&\frac{C}{\sqrt\e}\mE\int_0^T\Big|\bar{h}_1(X_t^{\e}, \sL_{X_t^{\e}})-\bar{h}_1(\bar{X}_t, \sL_{\bar{X}_t})-\sqrt\e\p_x\bar{h}_1(X_t^{\e}, \sL_{X_t^{\e}})U^\e_t\\
&&\qquad\qquad\qquad\qquad\qquad-\sqrt\e\tilde{\mE}[\p_\mu\bar{h}_1(X_t^{\e}, \sL_{X_t^{\e}})(\tilde{X}_t^{\e})\tilde{U}_t^{\e}]\Big|\dif t\\
&\leq&C\mE\int_0^T\big|[\p_x\bar{h}_1(\bar{X}_t+r\sqrt\e U^\e_t, \sL_{X_t^{\e}})-\p_x\bar{h}_1(X_t^{\e}, \sL_{X_t^{\e}})]U^\e_t\big|\dif t\\
&&+C\mE\int_0^T\left|\tilde{\mE}\Big[\Big(\p_\mu\bar{h}_1(\bar{X}_t, \sL_{\bar{X}_t+r\sqrt\e U^\e_t})(\tilde{\bar{X}}_t+r\sqrt\e \tilde{U}^\e_t)-\p_\mu\bar{h}_1(X_t^{\e}, \sL_{X_t^{\e}})
(\tilde{X}_t^{\e})\Big)\tilde{U}^\e_t\Big]\right|\dif t\\
&\leq&C\sqrt\e\int_0^T\mE|U^\e_t|^2\dif t\leq C_T\sqrt\e.
\de

For $J_3^\e$, following a similar calculation as for $J_2^\e$, we conclude that
\ce
J_3^\e&\leq&\frac{C}{\e}\mE\int_0^T\Big[\|\g_1(X_t^{\e}, \sL_{X_t^{\e}})-\g_1(\bar{X}_t, \sL_{\bar{X}_t})\|+\big\|\sqrt\e\p_x\g_1(X_t^{\e}, \sL_{X_t^{\e}})U^\e_t\\
&&\qquad\qquad\qquad+\sqrt\e\tilde{\mE}[\p_\mu\g_1(X_t^{\e}, \sL_{X_t^{\e}})(\tilde{X}_t^{\e})\tilde{U}_t^{\e}]\big\|\Big]\big\|\g_1(X_t^{\e}, \sL_{X_t^{\e}})-\g_1(\bar{X}_t, \sL_{\bar{X}_t})\\
&&\qquad\qquad\qquad-\sqrt\e\p_x\g_1(X_t^{\e}, \sL_{X_t^{\e}})U^\e_t-\sqrt\e\tilde{\mE}[\p_\mu\g_1(X_t^{\e}, \sL_{X_t^{\e}})(\tilde{X}_t^{\e})\tilde{U}_t^{\e}]\big\|\dif t\\
&\leq&C\sqrt\e\int_0^T\mE|U^\e_t|^3\dif t\leq C_T\sqrt\e.
\de

From the above deduction, it follows that
\ce
\Theta(\e)\leq C_T\sqrt\e.
\de
Thus the proof is complete.

\section{An example}\label{example}

Now let us give an example to illustrate the applicability of our result.
\bx\label{ex}
Consider the following slow-fast system:
\be\left\{\begin{array}{l}
\dif X_t^\e=\Big[\sin(aX_t^\e)+\int_{\mR}\check{x}\sL_{X_t^\e}(\dif\check{x})+\cos(bY_t^{\e, y_0, \sL_{\xi}})+\int_{\mR}\cos(q\check{y})\sL_{Y_t^{\e, \xi}}(\dif\check{y}) \Big]\dif t\\
\qquad\qquad+\int_{\mR}\sin(X_t^\e+\check{x})\sL_{X_t^\e}(\dif\check{x}) \dif B_t, \\
X_0^\e=\varrho, \quad 0\leq t\leq T, \\
\dif Y_t^{\e, \xi}=\frac{1}\e \int_{\mR}(-kY_t^{\e, \xi}+\l\check{y})\sL_{Y_t^{\e, \xi}}(\dif\check{y}) \dif t+\frac{1}{\sqrt\e}\int_{\mR}(mY_t^{\e, \xi}+\t\arctan\check{x})\sL_{X_t^\e}(\dif\check{x})\dif W_t, \\
Y_0^{\e, \xi}=\xi, \quad 0\leq t\leq T, \\
\dif Y_t^{\e, y_0, \sL_{\xi}}=\frac{1}\e\int_{\mR}(-kY_t^{\e, y_0, \sL_{\xi}}+\l\check{y})\sL_{Y_t^{\e, \xi}}(\dif\check{y}) \dif t+\frac{1}{\sqrt\e}\int_{\mR}(mY_t^{\e, y_0, \sL_{\xi}}+\t\arctan\check{x})\sL_{X_t^\e}(\dif\check{x})\dif W_t, \\
Y_0^{\e, y_0, \sL_{\xi}}=y_0, \quad 0\leq t\leq T,
\end{array}\right.
\label{exorieq}
\ee
where $\left(B_t\right),\left(W_t\right)$ are $1$-dimensional standard Brownian motions, respectively, defined on the complete filtered probability space $(\Omega,\sF,\{\sF_t\}_{t \in[0, T]}, \mathbb{P})$ and are mutually independent. $\varrho, \xi$ are $\sF_0$-measurable Gaussian random variables and $a,b,q,k,\l,m,\t$ are positive constants.
Let
$$h_1(x, \mu, y, \nu)=\sin(ax)+\int_{\mR}\check{x}\mu(\dif\check{x})+\cos(by)+\int_{\mR}\cos(q\check{y})\nu(\dif\check{y}),$$
$$\g_1(x,\mu)=\int_{\mR}\sin(x+\check{x})\mu(\dif\check{x}),$$
and
$$h_2(\mu, y, \nu)=-ky+\l\int_{\mR}\check{y}\nu(\dif\check{y}),\quad\g_2(\mu, y, \nu)=my+\t\int_{\mR}\arctan\check{x}\mu(\dif\check{x}).$$
For $x_i \in \mR, \mu_i \in \cP_2\left(\mR\right), i=1,2$, choosing $\pi^*\in\mathscr{C}(\mu_1,\mu_2)$ such that
\ce
&&|\g_1(x_1, \mu_1)-\g_1(x_2, \mu_2)|^2\no\\
&=&\left|\int_{\mR}\sin(x_1+\check{x}_1)\mu_1(\dif\check{x}_1)-\int_{\mR}\sin(x_2+\check{x}_2)\mu_2(\dif\check{x}_2)\right|^2\no\\
&\leq&\int_{\mR\times\mR}|\sin(x_1+\check{x}_1)-\sin(x_2+\check{x}_2)|^2\pi^*(\dif\check{x}_1,\dif\check{x}_2)\no\\
&\leq&2|x_1-x_2|^2+2\int_{\mR\times\mR}|\check{x}_1-\check{x}_2|^2\pi^*(\dif\check{x}_1,\dif\check{x}_2).
\de
By taking the infimum over $\sC(\mu_1,\mu_2)$, and using the characterization of the $L^2$-Wasserstein metric, we deduce that
\ce
|\g_1(x_1, \mu_1)-\g_1(x_2, \mu_2)|^2\leq2|x_1-x_2|^2+2\mW_2^2(\mu_1, \mu_2).
\de
By straightforward estimates one checks that $(\mathbf{H}_{h_1, \g_1}^1)$ holds. Moreover, for $\mu_i \in \cP_2\left(\mR\right)$, $y_i \in \mR$, $\nu_i \in \cP_2(\mR), i=1,2$,
\ce
&&|h_2(\mu_1, y_1, \nu_1)-h_2(\mu_2, y_2, \nu_2)|^2+|\g_2(\mu_1, y_1, \nu_1)-\g_2(\mu_2, y_2, \nu_2)|^2\\
&\leq&\max\{2k^2+2m^2,2\l^2,2\t^2\}\cdot(\mW_2^2(\mu_1, \mu_2)+|y_1-y_2|^2+\mW_2^2(\nu_1, \nu_2)),
\de
and
\ce
&&2\left\la y_1-y_2, h_2(\mu, y_1, \nu_1)-h_2(\mu, y_2, \nu_2)\right\ra+(2p-1)|\g_2(\mu, y_1, \nu_1)-\g_2(\mu, y_2, \nu_2)|^2\\
&=&2\left\la y_1-y_2, -ky_1+\l\int_{\mR}\check{y}_1\nu_1(\dif\check{y}_1)+ky_2-\l\int_{\mR}\check{y}_2\nu_2(\dif\check{y}_2)\right\ra+(2p-1)m^2|y_1-y_2|^2\\
&\leq&-2k|y_1-y_2|^2+2\l|y_1-y_2|\left|\int_{\mR}\check{y}_1\nu_1(\dif\check{y}_1)-\int_{\mR}\check{y}_2\nu_2(\dif\check{y}_2)\right|+(2p-1)m^2|y_1-y_2|^2\\
&\leq&-\left(2k-(2p-1)m^2-\l\right)|y_1-y_2|^2+\l\left|\int_{\mR}\check{y}_1\nu_1(\dif\check{y}_1)-\int_{\mR}\check{y}_2\nu_2(\dif\check{y}_2)\right|^2,
\de
which yields that
\ce
&&2\left\la y_1-y_2, h_2(\mu, y_1, \nu_1)-h_2(\mu, y_2, \nu_2)\right\ra+(2p-1)|\g_2(\mu, y_1, \nu_1)-\g_2(\mu, y_2, \nu_2)|^2\\
&\leq&-\left(2k-(2p-1)m^2-\l\right)|y_1-y_2|^2+\l\mW_2^2(\nu_1, \nu_2).
\de
If we choose appropriate $k$, $\l$, $m$ and $\t$, $(\mathbf{H}^1_{h_{2}, \g_{2}})$ and $(\mathbf{H}^2_{h_{2}, \g_{2}})$ are satisfied. For example, take $k=\frac{1}{8p}$, $\l=\frac{1}{32p}$, $m=\frac{1}{64p}$ and $\t=\frac{1}{64p}$. Note that
$$\p_\mu h_1(x, \mu, y, \nu)(\tilde{x})=1,\quad\p_\nu h_1(x, \mu, y, \nu)(\tilde{y})=-q\sin(q\tilde{y}),$$
$$\p_\mu\g_1(x,\mu)(\tilde{x})=\cos(x+\tilde{x}),\quad\p_\mu h_2(\mu,y,\nu)(\tilde{x})=0,\quad\p_\nu h_2(\mu,y,\nu)(\tilde{y})=\l,$$
$$\p_\mu\g_2(\mu,y,\nu)(\tilde{x})=\frac{\t}{1+\tilde{x}^2},\quad\p_{\tilde{x}}\p_\mu\g_2(\mu,y,\nu)(\tilde{x})=-\frac{2\t\tilde{x}}{(1+\tilde{x}^2)^2},$$
$$\p_{\tilde{x}}^2\p_\mu\g_2(\mu,y,\nu)(\tilde{x})=\frac{2\t(3\tilde{x}^2-1)}{(1+\tilde{x}^2)^3},\quad\p_{\tilde{x}}^3\p_\mu\g_2(\mu,y,\nu)(\tilde{x})=\frac{24\t\tilde{x}(1-\tilde{x}^2)}{(1+\tilde{x}^2)^4}.$$
A straightforward verification shows that $\g_1\in \big(C_b^{4,(1,3)}\cap C_b^{4,(2,2)}\cap C_b^{4,(3,1)}\big)(\mR\times\cP_2(\mR),\mR)$, $h_2, \g_2\in \big(C_b^{(3,1),6,(3,3)}\cap\mC_b^{(1,3),4,(2,2)}\cap\mC_b^{(2,2),2,(1,1)}\big)(\cP_2(\mR)
\times\mR\times\cP_2(\mR),\mR)$ and $h_1\in \big(C_b^{4,(3,1),6,(3,3)}\cap \mC_b^{4,(1,3),4,(2,2)}\cap\mC_b^{4,(2,2),2,(1,1)}\big)(\mR\times\cP_2(\mR)
\times\mR\times\cP_2(\mR),\mR)$. Consequently, the conclusions of Theorems \ref{xbarxp} and \ref{cltth} follow.
\ex


\begin{thebibliography}{999}
\bibitem{br}
V. Barbu and M. R\"{o}ckner: From nonlinear Fokker-Planck equations to solutions of distribution dependent SDE, {\it Ann. Probab.}, 48(2020)1902-1920.

\bibitem{br1}
V. Barbu and M. R\"{o}ckner: Uniqueness for nonlinear Fokker-Planck equations and weak uniqueness for McKean-Vlasov SDEs, {\it Stoch PDE: Anal Comp}, 9(2021)702-713.

\bibitem{bs}
Z. W. Bezemek and K. Spiliopoulos: Rate of homogenization for fully-coupled McKean-Vlasov SDEs, {\it Stoch. Dyn.}, 23(2023)2350013.

\bibitem{blpr}
R. Buckdahn, J. Li, S. Peng and C. Rainer: Mean-field stochastic differential equations
and associated PDEs, {\it Ann. Probab.}, 45(2017)824-878.

\bibitem{card}
P. Cardaliaguet: Notes on mean field games, P.-L. Lions Lectures at College de France, https://www.researchgate.net/publication/228702832. Notes on Mean Field Games.

\bibitem{carm}
R. Carmona and F. Delarue: Probability Theory of Mean Field Games with Applications I, Vol. 83. Springer, Cham, 2018.

\bibitem{ce}
S. Cerrai: Normal deviations from the averaged motion for some reaction-diffusion equations with fast oscillating perturbation, {\it J. Math. Pures Appl.}, 91(2009)614-647.

\bibitem{dq1}
X. Ding and H. Qiao: Euler-Maruyama approximations for stochastic McKean-Vlasov equations
with non-Lipschitz coefficients, {\it J. Theoret. Probab.}, 34(2021)1408-1425.

\bibitem{dq2}
X. Ding and H. Qiao: Stability for stochastic McKean-Vlasov equations with non-Lipschitz coefficients, {\it SIAM J. Control Optim.}, 59(2021)887-905.

\bibitem{fp}
J.P. Fouque, G. Papanicolaou, R. Sircar and K. Solna: Multiscale stochastic volatility for equity, interest rate, and credit derivatives, Cambridge University Press, Cambridge, 2011.

\bibitem{hlls}
W. Hong, S. Li, W. Liu and X. Sun: Central limit type theorem and large deviation principle for multi-scale McKean-Vlasov SDEs, {\it Probab. Theory Related Fields}, 187(2023)133-201.

\bibitem{hw}
X. Huang and F.-Y. Wang: Distribution dependent SDEs with singular coefficients, {\it Stochastic Process. Appl.}, 129(2019)4747-4770.

\bibitem{rk1}
R. Khasminskii: On stochastic processes defined by differential equations with a small parameter,
{\it Theory Probab. Appl.}, 11(1966)211-228.

\bibitem{kn}
C. Kumar and Neelima: On explicit Milstein-type scheme for McKean-Vlasov stochastic differential equations with super-linear drift coefficient, {\it Electron. J. Probab.}, 26(2021)1-32.

\bibitem{lwx}
Y. Li, F. Wu and L. Xie:
Poisson equation on Wasserstein space and diffusion approximations for multiscale McKean-Vlasov equation, {\it SIAM J. Math. Anal.}, 56(2024)1495-1524.

\bibitem{lx}
Y. Li and L. Xie: Functional law of large numbers and central limit theorem for slow-fast McKean-Vlasov equations, {\it Discrete Contin. Dyn. Syst. Ser. S}, 16(2023)846-877.

\bibitem{mte}
A. Majda, I. Timofeyev and E.V. Eijnden: A mathematical framework for stochastic climate models, {\it Comm. Pure Appl. Math.}, 54(2001)891-974.

\bibitem{hpm}
H.P. McKean: A class of Markov processes associated with nonlinear parabolic equations, Proc. Natl. Sci. U.S.A., 56(1966)1907-1911.

\bibitem{qw}
H. Qiao and W. Wei: Strong approximation of nonlinear filtering for multiscale McKean-Vlasov stochastic systems, https://arxiv.org/abs/2206.05037.

\bibitem{rrw}
P. Ren, M. R\"{o}ckner and F.-Y. Wang: Linearization of nonlinear Fokker-Planck equations and applications, {\it J. Differential Equations}, 322(2022)1-37.


\bibitem{rw}
P. Ren and F.-Y. Wang: Space-Distribution PDEs for path independent additive functionals of McKean-Vlasov SDEs, {\it Infin. Dimens. Anal. Quantum Probab. Relat. Top.}, 23(2020)2050018.

\bibitem{rsx}
 M. R\"{o}ckner, X. Sun and Y. Xie: Strong convergence order for slow-fast McKean-Vlasov stochastic differential equations, {\it Ann. Inst. H. Poincar\'{e} Probab. Statist.}, 57(2021)547-576.

\bibitem{rx}
M. R\"{o}ckner and L. Xie: Averaging principle and normal deviations for multiscale stochastic systems, {\it Commun. Math. Phys.}, 383(2021)1889-1937.

\bibitem{rxy}
M. R\"{o}ckner, L. Xie and L. Yang:
Asymptotic behavior of multiscale stochastic partial differential equations with H\"{o}lder coefficients, {\it J. Funct. Anal.}, 285(2023)110103.

\bibitem{rz}
M. R\"{o}ckner and X. Zhang: Well-posedness of distribution dependent SDEs with singular drifts, {\it Bernoulli}, 27(2021)1131-1158.

\bibitem{Wang}
F.-Y. Wang: Distribution dependent SDEs for Landau type equations, {\it Stochastic Process. Appl.}, 128(2018)595-621.

\bibitem{wr}
W. Wang and A.J. Roberts: Average and deviation for slow-fast stochastic partial differential equations, {\it J. Differential Equations}, 253(2012)1265-1286.

\bibitem{wtry}
F. Wu, T. Tian, J.B. Rawlings and G. Yin: Approximate method for stochastic chemical kinetics with two-time scales by chemical Langevin equations, {\it J. Chem. Phys.}, 144(2016)174112.

\bibitem{xq}
J. Xiang and H. Qiao: Averaging principles and central limit theorems for multiscale McKean-Vlasov stochastic systems, https://arxiv.org/abs/2501.11853.

\bibitem{xie}
L. Xie: Fast-slow stochastic dynamical system with singular coefficients, {\it Sci. China Math.}, 66(2023)819-838.

\end{thebibliography}
\end{document}